\documentclass[preprint,12pt]{elsarticle} 
\usepackage{b}
\usepackage[a4paper]{geometry}
\usepackage[small,it]{caption}
\usepackage{empheq}
\usepackage{tikz}
\usepackage{xifthen}
\usepackage{caption}
\usepackage{subfigure}
\usepackage{booktabs}

\usepackage{amsthm}

\newcommand{\bigo}{\mathcal{O}}
\newcommand{\x}{\mathbf{x}}
\newcommand{\y}{\mathbf{y}}
\newtheorem{remark}{Remark}

\newcommand*{\figfarnear}[1][1]{

\begin{tikzpicture}[scale=#1]

\foreach \x in {0,2}
{
\foreach \y in {0,2,4,6}
\draw (\x,\y) rectangle (\x + 2, \y + 2) node[midway] {$F$};
}

\foreach \x in {4,6}
{
\foreach \y in {0,2}
\draw (\x,\y) rectangle (\x + 2, \y + 2) node[midway] {$F$};
}

\draw (4,4) rectangle (5,5) node[midway] {$N$};
\draw (5,4) rectangle (6,5) node[midway] {$N$};
\draw (4,5) rectangle (5,6) node[midway] {$N$};
\draw (5,5) rectangle (6,6) node[midway] {$B$};
\draw (4,6) rectangle (5,7) node[midway] {$N$};
\draw (5,6) rectangle (6,7) node[midway] {$N$};
\draw (4,7) rectangle (5,8) node[midway] {$F$};
\draw (5,7) rectangle (6,8) node[midway] {$F$};
\draw (6,4) rectangle (7,5) node[midway] {$N$};
\draw (7,4) rectangle (8,5) node[midway] {$F$};
\draw (6,5) rectangle (7,6) node[midway] {$N$};
\draw (7,5) rectangle (8,6) node[midway] {$F$};

\draw (6,6) rectangle (8,8) node[midway] {$N$};

\end{tikzpicture}

}

\newcommand*{\figfarint}[1][1]{

\begin{tikzpicture}[scale=#1]

\foreach \x in {0,2,4,6}
{
\foreach \y in {0,2,4,6}
\ifthenelse{\x=0 \OR \y=0}
{\draw (\x,\y) rectangle (\x + 2, \y + 2) node[midway] {*};}
{\draw (\x,\y) rectangle (\x + 2, \y + 2);};
}
\draw (4,4) rectangle (5,5) node[midway] {};
\draw (5,4) rectangle (6,5) node[midway] {};
\draw (4,5) rectangle (5,6) node[midway] {};
\draw (5,5) rectangle (6,6) node[midway] {$B$};
\draw (4,6) rectangle (5,7) node[midway] {};
\draw (5,6) rectangle (6,7) node[midway] {};
\draw (4,7) rectangle (5,8) node[midway] {*};
\draw (5,7) rectangle (6,8) node[midway] {*};
\draw (6,4) rectangle (7,5) node[midway] {};
\draw (7,4) rectangle (8,5) node[midway] {*};
\draw (6,5) rectangle (7,6) node[midway] {};
\draw (7,5) rectangle (8,6) node[midway] {*};

\draw (6,6) rectangle (8,8) node[midway] {};

\end{tikzpicture}

}

\journal{Journal of Computational Physics}

\begin{document}

\begin{frontmatter}

\title{An adaptive fast multipole accelerated Poisson solver for complex geometries}
\author{T. Askham}
\author{A.J. Cerfon}
\address{Courant Institute of Mathematical Sciences, New York University, New York, NY 10012}
\date{\today}

\begin{abstract}
We present a fast, direct and adaptive Poisson solver for complex
two-dimensional geometries based on potential theory and fast multipole
acceleration. More precisely, the solver relies on the standard decomposition
of the solution as the sum of a volume integral to account for the source
distribution and a layer potential to enforce the desired boundary condition.
The volume integral is computed by applying the FMM on a square box that
encloses the domain of interest. For the sake of efficiency and convergence
acceleration, we first extend the source distribution (the right-hand side
in the Poisson equation) to the enclosing box as a $C^0$ function using
a fast, boundary integral-based method.  We demonstrate on
multiply connected domains with irregular boundaries that this continuous
extension leads to high accuracy without excessive adaptive refinement near
the boundary and, as a result, to an extremely efficient ``black box" fast
solver.
\end{abstract}


\end{frontmatter}

%%%%%%%%%%%%%%%%%%%%%%%%%%%%%%%%%%%%%%%%%%%%%%%%%%%%%%%%%%%%

\section{Introduction}
\label{sec:introduction} 
The solution of the Poisson equation is a critical task in many areas of
computational physics. The corresponding solvers need to be able to handle
complex, multiply connected geometries, to be fast, adaptive, and to yield
high order accuracy. Speed is of particular importance when the Poisson
equation is part of a larger system of equations or in the inner loop of an
optimization process. And since the physical quantity of interest is
often the gradient of the solution, rather than the
solution itself \cite{Gerris, Strauss, Celestin, Pataki, Lee}, partial derivatives of the solution must be computable with high accuracy as well.

Integral equation techniques have the potential to address all the challenges mentioned above. Complex geometries may be handled by decomposing the solution to Poisson's equation as the sum of a particular solution $v$ that does not satisfy the proper boundary condition in general, plus a homogeneous solution $u^H$ that solves Laplace's equation and is chosen so that the full solution $u=v+u^H$ satisfies the proper boundary condition. Fast and accurate solvers can be designed based on this construction. Indeed, several efficient and accurate integral equation based schemes exist to compute the solution of Laplace's equation on complex geometries \cite{Helsing08,Barnett1,Barnett2}, and fast and accurate evaluation of the particular solution $v$ on fully adaptive grids by use of the Fast Multipole Method (FMM) has also been demonstrated for Poisson's equation \cite{Ethridge, Langston}. Furthermore, in integral equation formulations derivatives do not have to be computed through direct numerical differentiation. Instead, one can analytically differentiate the kernels in the integral representation of the solution, and thus obtain integral representations for the derivatives of the solution as well. As a result the numerical error for the derivatives often converges at the same rate as the error for the solution itself \cite{Pataki, Langston}.
 
Remarkably, despite all the strengths described above, we are not aware of an integral equation based Poisson solver for planar problems that combines all the features at once. In \cite{Ethridge,JuneYub}, grid adaptivity and FMM acceleration are demonstrated, but only simple geometries are considered. In contrast, in \citep{McKenney}, Poisson's equation is solved for complex geometries and with FMM-accelerated quadratures, but the solver relies on fast methods for uniform grids \cite{Mayo1,Mayo2}. The purpose of this manuscript is to close this gap and to present an adaptive, FMM-accelerated Poisson solver for complex geometries. We achieve this in the following way. We embed the irregular domain $\Omega$ on which Poisson's equation needs to be solved in a larger square domain $\Omega_{B}$. We decompose the solution to Poisson's equation as $u=v+u^H$, and compute the particular solution $v$ on $\Omega_{B}$ with a fast and accurate solver for square domains \cite{Ethridge}. In order to calculate $v$ in this way, we need to extend the source function $f$ on the right-hand side of Poisson's equation beyond the domain $\Omega$ where it is originally given. We show that global function extension for $f$, constructed by solving Laplace's equation or a higher order partial differential equation on the domain $\mathbb{R}^2\setminus\Omega$, leads to a robust, efficient and accurate algorithm for the evaluation of $v$. This idea is very similar in spirit to the extension technique recently presented by Stein \textit{et al.}\cite{Stein2016} for the immersed boundary method, but quite different in its implementation. Our approach for computing $u^H$ is standard in its formulation \cite{Guenther88}, but it relies on numerical tools developed recently for optimized performance. Specifically, we represent $u^H$ as a layer potential whose density solves a second-kind integral equation. We use generalized Gaussian quadrature \cite{Bremer10,Bremer10U} to approximate the integrals, a fast direct solver \cite{Ho12} to compute the density and an FMM accelerated Quadrature By Expansion (QBX) algorithm \cite{Kloeckner13} to evaluate $u^H$ inside $\Omega$.

The structure of the article is as follows. In Section \ref{sec:potentialtheory} we present our formulation for the solution to Poisson's equation, which is based on standard potential theory. We stress its computational challenges, which are then addressed in the following sections. In section \ref{sec:boxcodes}, we describe an efficient and accurate algorithm for the evaluation of the particular solution $v$ and its derivatives in a square box. While this algorithm plays a central role in our approach to the problem, the section is relatively brief because our solver relies on an implementation of the algorithm and techniques that have been discussed in detail elsewhere \cite{Ethridge,Langston}. In section \ref{sec:irregular_box}, we explain how we use a global function extension algorithm in combination with a box Poisson solver for the computation of the particular solution $v$ on the whole square domain $\Omega_{B}$. This is a key element of our solver, which allows us to deal with complex geometries in an efficient manner. In Section \ref{sec:layer}, we present our numerical method for calculating the homogeneous solution $u^H$, as well as the function extension. Both are expressed as layer potentials and are computed in very similar ways. In Section \ref{sec:results} we study the performance of our new solver for two Poisson problems on a multiply connected domain. We summarize our work in Section \ref{sec:conclusion} and suggest directions for future work.

\section{The Potential Theoretic Approach to Poisson's Equation} 
\label{sec:potentialtheory}

In this article, we consider the solution $u$ to Poisson's equation with Dirichlet boundary conditions given by
\begin{align}
\Delta u &= f \mbox{ in } \Omega \label{eq:poisson1}\\
u &= g \mbox{ on } \partial \Omega \label{eq:poisson2}
\end{align}
where $\Omega$ is a smooth planar domain, which may or may not be multiply connected. The standard potential theory-based approach to the solution of (\ref{eq:poisson1}--
\ref{eq:poisson2}) proceeds as follows. 
The first step is to calculate a particular solution, i.e. a function $v$ which satisfies only equation (\ref{eq:poisson1}) but does in general not satisfy equation (\ref{eq:poisson2}). A natural candidate for $v$ is given by

\begin{equation}
\label{eq:volumeint1}
  v(\x) = \int_{\Omega} G(\x,\y) f(\y) \, d\y,
\end{equation}
where $G(\x,\y)$ is the free-space Green's function for Poisson's
equation. For planar problems, $G(\x,\y)=\log(||\x-\y||)/2\pi$. This is the situation we will consider in this article. Once $v$ has been computed, the second step is to compute a homogeneous solution with appropriate boundary conditions. Specifically, one solves the following Dirichlet problem

\begin{align}
\Delta u^H &= 0 \mbox{ in } \Omega \label{eq:laplace1}\\
u^H &= g-v|_{\partial \Omega} \mbox{ on } \partial \Omega \label{eq:laplace2}.
\end{align}
The solution to (\ref{eq:poisson1}--\ref{eq:poisson2}) is then the sum,
$u = v+u^H$. There are many options for the numerical implementation of these two steps and we will not attempt to provide an exhaustive review of them here. Instead, we focus on our new approach, which is designed to address situations for which the domain $\Omega$ may be irregular and where derivatives of the solution are also required with high accuracy. The purpose of this section is to give a short overview of our approach. This overview is divided into two subsections: subsection \ref{sec:part_summary} concerns the computation of $v$, and subsection \ref{sec:hom_summary} describes the computation of $u^H$. The presentation in these two subsections is meant to give a general idea of our numerical scheme, and is brief on purpose. We provide detailed descriptions of our numerical methods to calculate of $v$ and $u^H$ in sections \ref{sec:boxcodes} and \ref{sec:irregular_box} for $v$, and section \ref{sec:layer} for $u^H$.

\subsection{Computing the particular solution}
\label{sec:part_summary}
There are two challenges associated with the evaluation of the particular solution $v$ through the integral (\ref{eq:volumeint1}).
First, accurate quadratures must be used in order to handle the logarithmic singularity. Second, given a quadrature rule, the na\"{i}ve numerical approach to computing (\ref{eq:volumeint1}) would require $\bigo (N^2)$ work for a domain with $N$ discretization nodes. It is now well known that the computational work can in fact be reduced to $\bigo (N)$ via the fast multipole method \cite{Carrier88}. Furthermore, for a fixed domain $\Omega$, the quadrature rules for  a weakly singular kernel $G(\x,\y)$ can be precomputed using an adaptive, brute-force procedure \cite{Genz80}. As a result, there exist particularly efficient $\bigo (N)$ algorithms, including optimized versions of the FMM \cite{Ethridge}, to compute the integral (\ref{eq:volumeint1}) for problems specified on a box. We choose such an algorithm for our solver, and provide some details of this type of method in Section \ref{sec:boxcodes}.

%Furthermore, for a fixed, regular domain $\Omega$ (e.g. a square or disk), the construction of accurate quadrature rules for the weakly singular kernel  $G(x,y)$ is well understood, see, \em{inter alia}, \cite{Genz80}\color{red}REFERENCES? \color{black}. 

For irregular domains $\Omega$, however, the situation is quite different. The calculation of appropriate quadratures is much more difficult and fewer optimizations of the FMM are available. A natural strategy, then, for irregular domains is to consider a larger, square domain $\Omega_B$ containing $\Omega$ and to instead compute

\begin{equation}
\label{eq:volumeint2}
  v(\x) = \int_{\Omega_B} G(\x,\y) f_e(\y) \, d\y,
\end{equation}
where $f_{e}$ is defined on all of $\Omega_B$, and constructed such that $f_e = f$ on $\Omega$. One of the main novelties of our work is to compute $f_e$ via global function extension: $f_{e}$ restricted to $\R^2\setminus\Omega$ is the solution of an elliptic partial differential equation with Dirichlet data $f_{e}=f$ on $\partial\Omega$.The PDE is solved with a standard integral equation representation. We elaborate on this idea in Section \ref{sec:irregular_box}.

\begin{remark}
It should be noted here that the solution provided by any Poisson solver is a valid particular solution, though it will not necessarily be equal to the one given by \eqref{eq:volumeint1}. This includes in particular the solutions produced by FFT-based solvers for rectangular and circular domains, which are very fast in terms of 
work per grid point. The method of \cite{McKenney} uses such a particular solution, computed via Buneman's method \cite{Buzbee70} and the modified stencils developed in \cite{Mayo1}.
The algorithm of \cite{Ethridge}, which we choose for our solver, is an  alternative to such methods, with its greatest advantage being  the ease with which it handles adaptive discretization. We demonstrate this advantage with numerical examples in Section \ref{sec:results}.
\end{remark}

\subsection{Computing the homogeneous solution}
\label{sec:hom_summary}

A standard approach to the solution of Laplace's equation
is to represent the solution $u^H$ as a layer potential with 
unknown density $\mu$ on the  boundary. The representation should be 
chosen so that imposing the boundary conditions results in an 
invertible, second kind integral equation (SKIE) for the density on the boundary.
This is a well-studied area and there exist appropriate integral representations
for multiply connected domains \cite{Greenbaum93,Helsing05}, unbounded domains \cite{Guenther88}, and for situations with other types of boundary conditions
\cite{Ojala12}. Further references can be found in the previously
cited papers, and we recommend \cite{Guenther88,Atkinson} for very clear treatments of this topic. 

Once a suitable representation for $u^H$ is chosen, the discretization of the problem is then simply a matter of quadrature for the resulting SKIE. In general, the integral kernel may be singular and the choice of quadrature requires attention \cite{Alpert99,Kapur,Kloeckner13,Helsing08,Bremer10,Bremer10U}.
Once discretized, there are many tools available for the fast solution of the resulting linear system, which we briefly discuss in Section \ref{sec:layer}. For this article, we choose a direct method \cite{Ho12} that is optimized for the type of problems considered here. 

After the density $\sigma$ is computed, the solution $u^H$ can
be evaluated in the domain. This step is  trivially direct but it is not without its difficulties. With $N$ discretization points in the domain and $M$
discretization points on the boundary, na\"{i}ve 
computation of the necessary integrals would require
$\bigo (MN)$ work. This work can be reduced to $\bigo
(M+N)$ with the FMM.
Because the integral kernel of the solution representation
is typically singular in the ambient
space ($\R^2$), computing the potential
$u^H$ to high accuracy near the boundary requires special
quadrature schemes. Such schemes have been developed recently \cite{Helsing08,Kloeckner13}, and for our solver we choose to rely on the Quadrature By Expansion method (QBX) \cite{Kloeckner13}, which we also briefly discuss 
in section \ref{sec:layer}.

\section{Box Codes}
\label{sec:boxcodes}

This section reviews relevant features of the algorithm of \cite{Ethridge}, which is the original ``box code'', and which we have implemented in our solver. By ``box codes'', we mean a class of fast solvers which are  used to evaluate integrals of the form

\begin{equation}
Vf(\x) = \int_\Omega G(\x,\y) f(\y) \, d\y \; , \label{eq:volumeint3} 
\end{equation}
where the integral kernel $G(\x,\y)$ is a translation 
invariant Green's function, the domain $\Omega$ is a box, 
and $f(\y)$ is a given density. We take 
$G(\x,\y) = -\log \|\x-\y\|/2\pi$ in what follows.

\subsection{Outline of a box code}
\label{sec:box_outline}

As in all fast multipole methods, a FMM-based box code is 
based on a hierarchical division of space.
Specifically, the domain is taken to be the root box 
(level $0$) of a quad-tree. The finer levels are obtained
by subdividing boxes from the previous level into 
four equal parts. After a box is subdivided, the four resulting
boxes on the next level are its children. The quad-tree
for a box code is thus fairly typical for an FMM.
The primary distinction of a quad-tree as used in a box
code is that it is typically a level restricted tree, i.e., 
adjacent leaf boxes are required to be no more than one 
level apart in the tree hierarchy. 

When computing \eqref{eq:volumeint3}, a choice has to be made
as to how the function $f$ is represented on each leaf
box of the quad-tree. The standard choice, as in
\cite{Ethridge}, is to represent $f$ by collocation
points (for monomials, Chebyshev polynomials, etc) on
each leaf box, using the same points scaled for each
level. Then, a reasonable subdivision criterion
for a box is whether or not the function $f$ is well
approximated by its interpolant up to a given tolerance
on that box. This criterion makes a box code an adaptive
method, with the order of accuracy determined by the order
of the polynomial approximation on each box.

After the quad-tree is formed, we have that 
$\Omega = \cup_j B_j$ where the $B_j$ are leaf boxes and 
on each leaf box there is a polynomial $p_j$ which
approximates the density $f$. Let $\tilde{f}$, defined
by setting $\tilde{f}(\x) = p_j(\x)$ 
for $\x \in B_j$, be the approximation of $f$ over the
whole domain.
The box code proceeds to evaluate the potential 
$V\tilde{f}(x)$, where the evaluation points $x$ are 
taken to be the collocation points 
of the polynomials $p_j$. Let $\tilde{V}\tilde{f}(\x)$ be the
computed values of $V\tilde{f}(\x)$. To evaluate the
volume integral at other points in the domain,
we evaluate the polynomial which interpolates the 
values $\tilde{V}\tilde{f}(\x)$ on each box. We denote
this piecewise polynomial function by $\tilde{v}$.
The distinction between $\tilde{v}$ and $\tilde{V}\tilde{f}$
is subtle but necessary here. For the sake of speed,
a box code only evaluates $\tilde{V}\tilde{f}$ at
the collocation nodes. The operator $V$ is approximated
more or less exactly so the error is determined by the
interpolation error for $\tilde{f}$.
The values of $\tilde{v}$ incur further interpolation
error, which depnds on the order of the nodes and the
smoothness of $\tilde{V}\tilde{f}$. We address the error
analysis in more detail in Section \ref{subsec:smoothf}.

For a quad-tree
with $N$ total collocation points, computing $\tilde{V}\tilde{f}$ 
would require $\bigo (N^2)$ operations if done 
na\"{i}vely. This cost can be reduced to $\bigo (N)$ by 
using the fast multipole method. In the context of this article, it
is only necessary to describe the result of the FMM.
For a detailed account of the structure of the FMM, see 
\cite{Carrier88,Hrycak98,Ethridge}. 

\begin{figure}[h!]
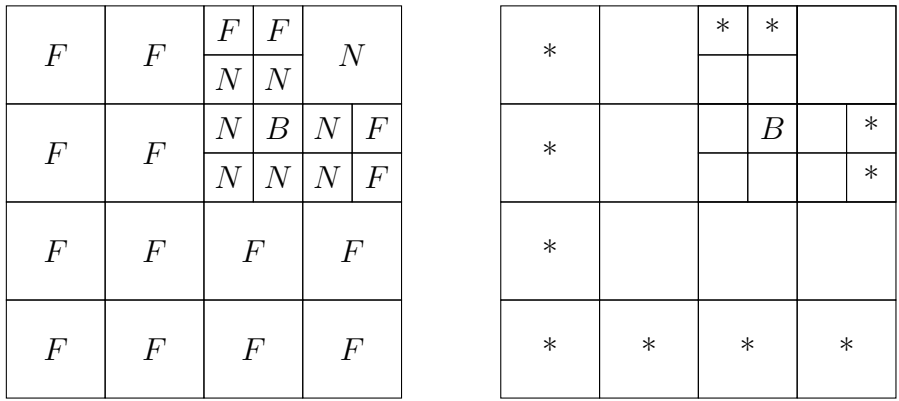
 

\centering

\begin{minipage}{0.4\textwidth}
\figfarnear[.65] 
\end{minipage}
\begin{minipage}{0.4\textwidth}
\figfarint[.65]
\end{minipage}

\caption{In the figure on the left, the leaves of a quad-tree
are shown and the boxes in the near field of the box $B$ 
are marked with an $N$ while the boxes in the far field of 
$B$ are marked with an $F$.
The same quad tree is shown on the right and the boxes for 
which $B$ is in the far field are marked with an asterisk (*),
these boxes being in $\mathcal{F}(B)$.}

\label{fig:farnear}

\end{figure}

Let $B_j$ be a leaf box of the
quad-tree with width $h$. The ``near field'' of $B_j$ is defined 
to be any leaf box whose interior intersects the interior of 
the box of width $3h$ centered at $B_j$. The boxes which are not 
in the near field of $B_j$ are said to be in the ``far field.''
Because the boxes in the far field of $B_j$ are separated 
from $B_j$ by a box of at least the same size as $B_j$, these
boxes are said to be ``well separated.''
Let 
$\mathcal{F}(B_j) = \{ i : B_j \mbox{ is in the far field of } B_i \}$
be the set of leaf boxes for which $B_j$ is well separated 
and $\Omega_j = \cup_{i \in \mathcal{F}(B_j)} B_i$ be the union over
these leaf boxes. For a non-uniform tree, it is not necessarily the
case that the boxes of $\mathcal{F}(B_j)$ and the far field of
$B_j$ are the same. See Figure \ref{fig:farnear} for examples of these
sets. In $\bigo (N)$
time, the FMM computes functions $\Phi_j$ 
for each leaf box $B_j$ which are expansions (more precisely, 
the sum of a Taylor expansion
and a number of multipole expansions) 
approximating the influence of all leaf boxes in $\mathcal{F}(B_j)$ at any 
point in $B_j$, i.e. for any $\x \in B_j$

$$ \Phi_j(\x) \approx \int_{\Omega_j} G(\x,\y) \tilde{f} (\y) \, d\y \, . $$
With $\Phi_j$ computed, it is possible to compute the
volume integral \eqref{eq:volumeint3} by directly adding the 
influence of leaf boxes for which $B_j$ is in the near field, 
i.e. for any $\x \in B_j$

$$ V \tilde{f} (\x) \approx \Phi_j(\x) + 
\sum_{i \not\in \mathcal{F}(B_j)} \int_{B_i} G(\x,\y) \tilde{f} (\y) \, d\y. $$
where the second term on the right-hand side 
is evaluated by direct computation, using a high order 
quadrature rule. This step is $\bigo (1)$ per point because the number of 
boxes for which $B_j$ is in the near field and the cost of evaluating
$\Phi_j$ are bounded independent
of $N$. For a given precision $\eps_V$, the computed values
$$ \tilde{V}\tilde{f}(\x) = \Phi_j(\x) + 
\sum_{i \not\in \mathcal{F}(B_j)} \int_{B_i} G(\x,\y) \tilde{f} (\y) \, d\y. $$
satisfy

$$ |\tilde{V}\tilde{f}(\x) - V\tilde{f}(\x) | \leq \eps_V \|\tilde{f}\|_{L^1}. $$
To achieve this bound for smaller values of $\eps_V$, the FMM 
uses higher-order expansions to approximate $\Phi_j$. See \cite{Greengard88} 
for more on the error analysis of the FMM. 

While $\bigo (N)$ is indeed optimal in terms of complexity,
the numerical scheme presented in \cite{Ethridge} is
particularly fast in terms of work per gridpoint. 
For far field interactions, the speed is due in part to 
the fact that the translation
of multipole expansions is diagonalized through the use of 
plane wave expansions,
see \cite{Ethridge} and \cite{Hrycak98} for details.
For near field interactions, the speed is due to the use 
of precomputed tables. Because the tree is level-restricted,
there are a limited number of near field interactions
possible, up to scale. Therefore, if the possible
interactions are stored for a unit box, the influence of 
any box on a box in its near field can be computed at the cost
of a small matrix-vector multiply.

\subsection{Derivatives of the potential}
\label{subsec:volderivatives}
In many physical applications, the derivatives of the volume 
potential $Vf(\x)$ are the quantities of interest, instead of $Vf(\x)$ itself. Once
the values of the potential $\tilde{V}\tilde{f}$ are computed,
one could differentiate the piecewise polynomial function, 
$\tilde{v}$, which interpolates
the potential on each leaf box to obtain an approximation of the 
derivatives. This computation results in derivative values 
which have an order of accuracy that is one less than the order
of accuracy for the potential. 

Instead, the derivatives can be computed by recognizing that 
they are given by another volume integral, e.g.

\begin{equation}
\partial_{x_1} Vf(\x) =
\int_\Omega \partial_{x_1} G(\x,\y) f(\y) \, d\y 
\, .
\label{eq:volumeint3x}
\end{equation}
In fact, the volume integral for the derivatives can be 
computed alongside the evaluation of the volume integral
for the potential with modest impact on the run time. As
in the case of computing the potential, the near-field interactions 
can be calculated making use of precomputed tables. The
far-field interactions can be computed by differentiating 
the local expansion 
for the far-field, i.e. by differentiating $\Phi_j(\x)$, 
which is typically
a much higher order approximation than the order of the 
collocation points. The result of computing the derivatives
of the solution with this approach is that the derivatives 
display the same convergence rate as the potential. For the
calculations presented in this article, the authors have 
implemented such a scheme. A similar approach to computing 
derivatives was taken in \cite{Langston}.

\subsection{Error analysis for smooth $f$}
\label{subsec:smoothf}
As above, let $\tilde{f}$ denote the piecewise
polynomial approximation to $f$ for a given
tree and let $\tilde{V}\tilde{f}$ denote the
computed value of $V\tilde{f}$.
Suppose that the fast multipole method is applied 
with precision $\eps_V$ and that the local 
interaction tables are computed to at least that 
precision. Then, the error in $V\tilde{f}$ 
at a collocation node $\x$ has the following bound:

\begin{equation}
  |\tilde{V}\tilde{f}(\x)-\tilde{V}\tilde{f}(\x)| 
  \leq \eps_V \|\tilde{f}\|_1 \; .
\end{equation}
That is, the values of $V\tilde{f}(\x)$ 
are computed at the collocation nodes with an error that 
depends only on the truncation order of the fast multipole
method. It then follows that the total error at any
given collocation node $\x$ is bounded by

\begin{equation}
 | Vf(\x) - \tilde{V}\tilde{f}(\x)| \leq 
|V\tilde{f}(\x)-\tilde{V}\tilde{f}(\x) |
+ |V(\tilde{f}-f)(\x)| \leq \eps_V \|\tilde{f}\|_1+ 
C_\Omega \|\tilde{f}-f\|_\infty 
\, , \label{eq:errorAnalysis}
\end{equation}
for a constant $C_\Omega$ which is independent of $f$ and 
given by

\begin{equation}
 C_\Omega = \max_{\x \in \Omega} \int_\Omega |G(\x,\y)| \, d\y \, .
\end{equation}
The bound \eqref{eq:errorAnalysis} provides an \textit{a priori}
estimate of the accuracy of the solution which depends only
on the values of the data $f$. This is useful when designing
adaptive refinement strategies as one can simply check whether
$f$ is well approximated on each leaf in the tree. On a uniform
tree, we see that the order of accuracy of the overall scheme
depends on the order of accuracy of the local polynomial approximation
to $f$ on leaf nodes. Finally, we note that \eqref{eq:errorAnalysis} 
only depends on the fact that $V$ is a bounded operator from 
$L^\infty$ to $L^\infty$. In particular, it is clear that 
\eqref{eq:errorAnalysis} holds analogously for $\nabla V$
and that on a uniform tree we should see the same order of
accuracy for the potential and gradient values (this is 
sometimes referred to as ``super-convergence'').

In some cases, the values of the potential will be desired 
at points other than the collocation nodes. These values 
can be obtained by interpolation of the computed potential
values, i.e. by evaluating $\tilde{v}$. Generally, the error 
in this step will be of the same
order as the error in \eqref{eq:errorAnalysis} because the
potential is a smoother function than $f$ as a function on
$\R^2$. However, that reasoning does not always apply when
$f$ is a smooth function on $\Omega$. This is because the 
volume integral \eqref{eq:volumeint3} implicitly defines 
$f$ to be zero outside of $\Omega$. 

By example, the function
$f \equiv 1$ on $\Omega = [-1,1] \times [-1,1]$ is
well-resolved by a tree with a single leaf box 
(equal to $\Omega$ itself) containing, say, a $4\times 4$ tensor
grid of equispaced points (though virtually any grid would
work here). The box code would produce the correct 
values of \eqref{eq:volumeint3} to near machine precision
at each collocation node. The problem is that the resulting
function is not well-resolved by the grid and the 
error in interpolating the function at arbitrary points 
in $\Omega$ is large. This simple example demonstrates an 
important point: the global smoothness of $f$ matters.

\section{Box codes for irregular domains}
\label{sec:irregular_box}

Perhaps the most natural idea to compute
a particular solution $v$ to Poisson's equation for an arbitrary 
irregular domain by using a box code is as follows. Suppose $\Omega$ 
is the irregular domain and $\Omega_B$ 
is a box such that $\Omega \subset \Omega_B$. Then, a particular 
solution on $\Omega$ can be computed as in \eqref{eq:volumeint2}, 
so long as an extension $f_e$ on $\Omega_B \setminus \Omega$ of the 
right hand side $f$ is given. 

In many cases of practical interest, a smooth or continuous
extension $f_e$ of $f$ is readily available. The density 
$f$ may for example describe a compactly supported distribution of electric charges which smoothly goes to zero on the boundary of $\Omega$. Likewise, the magnetohydrodynamic equilibrium of a plasma
confined in a tokamak is given by a semilinear Poisson equation in which the right-hand side $f$ smoothly goes to zero on the boundary of $\Omega$ in most situations \cite{Pataki,Lee}. In such cases, $f_{e}\equiv 0$ on $\Omega_B \setminus \Omega$ is a natural and satisfying choice.

In the general case, however, $f_{e}$ does not have an obvious physical meaning, and the extension $f_e$ is constructed as a purely mathematical artifice required by the box solver. The problem of specifying a function $f_e$ such that $f_e = f$ on $\Omega$ is extremely open; we narrow it by looking for an extension $f_e$ that is favorable in terms of the efficiency and accuracy of the box code.

\subsection{Error analysis for non-smooth densities}
\label{subsec:nonsmoothf}
To motivate our construction of $f_{e}$, we first perform a heuristic but more detailed analysis  of the error bound \eqref{eq:errorAnalysis} for densities $f_e$ which 
are not necessarily smooth on the box $\Omega_B$.
Using the standard multi-index notation, let 
$\partial^{\boldsymbol{\alpha}} = \partial_1^{\alpha_1} \partial_2^{\alpha_2}$
and $|\boldsymbol{\alpha}| = \alpha_1 + \alpha_2$. Then the
differentiability class $C^k(A)$ of a domain $A$ is defined 
to be the set of functions $g$ such that 
$\partial^{\boldsymbol{\alpha}} g$ is continous for each 
$\boldsymbol{\alpha}$ with $|\boldsymbol{\alpha}| \leq k$,
with the convention that $C^{-1}(A)$ is the set of 
bounded functions which are possibly discontinuous.
For a density $f_e \in C^{k}(\Omega_B)$ and
a uniform tree with leaf boxes of width $h$, let $\tilde{f_e}$ 
be the numerical approximation as in the previous section, using
interpolants of order $p$ (degree $p-1$) on each box. Standard
error estimates imply that 

\begin{equation}
\|f_e-\tilde{f_e}\|_\infty = \bigo (h^m) \, ,
\end{equation}
where $m = \min (k,p)$. 
If the additional assumption is made that $f_e$ is piecewise $C^{k+l}$
for some $l>0$ (say that for a domain $A \subset \Omega_B$, 
the density $f_e \in C^{k+l}(A)$ and $f_e \in C^{k+l}(\Omega_B 
\setminus A)$) then the approximation order is improved to 
$m = \min(k+1,p)$, see, {\em inter alia}, \cite{Dragomir08}.

These bounds suggest that the scheme should have
$\bigo (1)$ error for a piecewise smooth density which
is discontinuous across some boundary. However, in practice the observed
convergence rate is faster, even for the derivatives
of $Vf_e$. The reason for this is that the bound \eqref{eq:errorAnalysis} only makes
use of the fact that $V$ and its derivatives are bounded on $L^\infty$. It does not take into account the fact that they  are given by weakly singular integrals. One could thus seek a tighter bound for the error $|Vf_e(\x)-V\tilde{f_e}(\x)|$, but with our construction of the solution given by Equations \eqref{eq:volumeint1}, \eqref{eq:laplace1} and \eqref{eq:laplace2}, it suffices to see how good of a particular solution 
$V\tilde{f_e}$ is. For this purpose, let $\x$ be contained in a box 
$B_j$ and let $\Omega_j$ denote the union over all boxes for which $B_j$ is in the far-field, as already defined in Section \ref{sec:box_outline}. We can write

\begin{eqnarray}
  V\tilde{f_e}(\x) &=& \int_\Omega G(\x,\y) \tilde{f_e}(\y) \, d\y \\
  &=& \int_{\Omega\setminus \Omega_j} G(\x,\y) \tilde{f_e}(\y) \, d\y
  + \int_{\Omega_j} G(\x,\y) \tilde{f_e}(\y) \, d\y \, .
\end{eqnarray}
The contribution to $V\tilde{f_e}(\x)$ from the second term is harmonic. The contribution
from the first term is the relevant one regarding the
quality of $V\tilde{f_e}$ as a particular solution.
Let $h$ be the side length of the box $B_j$. If we consider 
the error without the far-field contribution, we have

\begin{eqnarray}
 \left | \int_{\Omega\setminus \Omega_j} G(\x,\y) \tilde{f_e}(\y) \, d\y
-\int_{\Omega\setminus \Omega_j} G(\x,\y) f_e(\y) \, d\y \right |
&\leq& \|\tilde{f_e}-f_e\|_\infty \int_{\Omega\setminus \Omega_j} |G(\x,\y)| \, d\y  \\
&=& \bigo (h^2 |\log h|) \|\tilde{f_e}-f_e\|_\infty \; . \label{eq:potBound}
\end{eqnarray}
For the derivative values we have the analogous bound

\begin{eqnarray}
 \left | \int_{\Omega\setminus \Omega_j} \nabla G(\x,\y) \tilde{f_e}(\y) \, d\y
-\int_{\Omega\setminus \Omega_j} \nabla G(\x,\y) f_e(\y) \, d\y \right |
&\leq& \|\tilde{f_e}-f_e\|_\infty \int_{\Omega\setminus \Omega_j} 
|\nabla G(\x,\y)| \, d\y  \\
&=& \bigo (h) \|\tilde{f_e}-f_e\|_\infty \; . \label{eq:gradBound}
\end{eqnarray}

There are two main conclusions from the preceding analysis. The first conclusion is the intuitive result that the smoother $f_e$ is the better the approximation of the particular solution. The more interesting conclusion is that $f_e$ may not need to be quite as smooth as we may have initially expected. Specifically, in our implementation, we use 4th order interpolants on each leaf box. Suppose we discretize the domain 
with a uniform tree and leaf boxes of side length $h$.
For a piecewise smooth density $f_e$ which is discontinuous, 
the above bounds imply (nearly) 2nd order accuracy in the values of
the particular solution and 1st order accuracy in the gradient. 
Similarly, for a piecewise smooth density $f_e$ which is continuous,
they imply (nearly) 3rd and 2nd order accuracy, respectively.
These bounds are consistent with our numerical results, as we demonstrate in Section \ref{sec:results}.

Finally, note that \eqref{eq:potBound} and \eqref{eq:gradBound} also imply that we do not expect to observe ``super-convergence'' for piecewise
smooth densities unless they are of sufficient smoothness
globally. While super-convergence would be a desirable
property, we have found that adaptive refinement strategies
can be advantageously used to obtain the desired high-accuracy (though not necessarily high-order accurate) values for derivatives of the potential. We present these numerical results in more detail in Section \ref{sec:results}.

\subsection{Global function extension}
\label{subsec:funextend}

In previous attempts to use box codes for irregular domains,
two main extrapolation techniques were used. The first,
which we call ``extrapolation by zero'', simply sets the
density $f_e$ to be zero outside the domain\cite{Ethridge}. In this case,
the function $f_e$ is as smooth as the original density $f$ 
inside the domain and trivially smooth outside the domain. Therefore, the estimates for piecewise smooth functions apply and we see that the scheme should converge
with a rate $\bigo (h^2|\log h|)$ for a uniform tree. The reader may keep in mind that for such situations, a box code relying on an adaptive tree is more efficient in terms of degrees of freedom than a code using a uniform tree. Even if so, adaptive refinement for  functions with a discontinuity requires trees with a large number of refinement levels and therefore a large number of grid points. This can make the ``extrapolation
by zero'' approach computationally costly. The second
extrapolation method uses local polynomial approximations to $f$
to extrapolate $f$ outside the domain over short distances \cite{Langston2}.
A major limitation of this method is that it results in a smooth $f_e$ for individual leaf boxes but has no guarantees of smoothness across boxes. Since there
can be discontinuities in $f_e$ across boxes near the original domain, the computed potential $Vf_e$ may be unresolved on those boxes. This issue seems to be inherent in local
extrapolation methods. That is why we seek out
a global extrapolation method which improves on the na\"{i}ve global ``extrapolation
by zero'' approach, in terms of both speed and accuracy. 

The global extrapolation method we adopt is similar in spirit to the one recently proposed in \cite{Stein2016} in a different context, and can be explained in a few words. Let $w$ solve the PDE

\begin{equation}
\begin{array}{rcl} 
  \Delta w = 0 & \mbox{ in } & \R^2 \setminus \Omega \\
  w = f & \mbox{ on } & \partial \Omega
\end{array}
\, , \label{eq:ctsExtrapPDE}
\end{equation}
subject to the condition that $w(\mathbf{x})$ is bounded as
$\|\mathbf{x}\| \to \infty$. Then, the function $f_e$
defined by

\begin{equation}
\begin{array}{rcl}
f_e(\mathbf{x}) = f(\mathbf{x}) & \mbox{ for } 
& \mathbf{x} \in \Omega \\
f_e(\mathbf{x}) = w(\mathbf{x}) & \mbox{ for } 
& \mathbf{x} \in \Omega_B 
\setminus \Omega
\end{array}
\end{equation}
is globally continuous, as smooth as $f$ on $\Omega$, and 
smooth on $\Omega_B \setminus \Omega$. While this may at first seem 
like a computationally expensive way to extrapolate $f$,
the analytical and numerical machinery required to solve this
problem is in fact the same as that required to solve the
harmonic problem \eqref{eq:laplace1}--\eqref{eq:laplace2},
which is used to enforce the boundary condition of the original
Poisson problem (\eqref{eq:poisson1}--\eqref{eq:poisson2}).

Note that the method can be generalized to compute globally
$C^k$ extrapolations of $f$ by solving polyharmonic 
equations. For example, a $C^1$ extrapolation can be 
computed by solving the following biharmonic problem:

\begin{equation}
\begin{array}{rcl} 
  \Delta^2 w = 0 & \mbox{ in } & \R^2 \setminus \Omega \\
  w = f & \mbox{ on } & \partial \Omega \\
  \dfrac{\partial w}{\partial n} = \dfrac{\partial f}{\partial n} 
& \mbox{ on } & \partial \Omega \\
  w = 0 & \mbox{ on } & \partial \Omega_e \\
  \dfrac{\partial w}{\partial n} = 0 
& \mbox{ on } & \partial \Omega_e 
\end{array}
\, , \label{eq:c1ExtrapPDE}
\end{equation}
where $\Omega_e$ is some domain containing $\Omega$. 
Once $w$ is computed, $f_e$ is as defined in the continuous case. 
While the methods for Equation \eqref{eq:c1ExtrapPDE} are not as well developed as in
the Laplace case, there exist similar potential-theory based
integral equations and fast solution methods for the solution 
of the biharmonic problem. See, for instance, P. Farkas' PhD thesis 
\cite{Farkas89} and the approach of \cite{Askham}
which leverages the Sherman-Lauricella integral equations for
elasticity. There are two main reasons to not consider 
extrapolations based on polyharmonic equations of higher order
than the biharmonic equation: (1) very few numerical tools have been developed for such equations
and (2) the equations require to provide high order 
derivatives of the data $f$ in the direction normal to the boundary, which in most physical applications are not readily available, and can be challenging to compute with high accuracy numerically, even when using integral equation based methods \cite{Ricketson}.

For our numerical tests, and in the version of the code which will be available online, only the harmonic expansion calculated by solving \eqref{eq:ctsExtrapPDE} is implemented. The details of our implementation are given in \ref{sec:layer} in parallel with the calculation of $u^H$, since both computations rely on the same mathematical and numerical tools.

\section{Computing the homogeneous solution and the harmonic extension}
\label{sec:layer}

In this section, we describe how we compute the homogeneous solution $u^{H}$ which solves the harmonic problem (\eqref{eq:laplace1}--\eqref{eq:laplace2}). Since we use similar numerical techniques to solve this problem and to compute the global function extension through \eqref{eq:ctsExtrapPDE}, we will also discuss the latter, and highlight the small differences between the two situations. 

\subsection{Layer potentials}
\label{sec:layer_representation}

Before we proceed, we should clarify what we mean by
a multiply connected domain and the normal direction
to the boundary curve. Let $\Omega$ be an interior
domain with boundary $\partial\Omega$. For a multiply connected
domain, $\partial\Omega$ is given as the union over disjoint, closed 
curves $\partial\Omega = \bigcup_{i=0}^l \Gamma_i$, with $\Gamma_0$
corresponding to the outer boundary. The normal direction on
each component $\Gamma_i$ is taken to be the direction pointing 
away from $\Omega$. For $\Gamma_0$, this vector points to the exterior of the
curve and for $\Gamma_i$, $i = 1,\ldots,l$, this points to the interior
of the curve. To see a simple illustration of such a domain and 
its normal vectors, see Figure \ref{fig:omega} in Section \ref{sec:results}.

We write both the homogeneous solution $u^H$ and the extension $w$ of $f$ in $\Omega_{B}\setminus\Omega$ as layer potentials \cite{Guenther88}. Specifically, for the homogeneous solution we write
\begin{equation}
 u^H (\x) = S \mu (\x) + D \mu (\x)
 \label{eq:layer_for_uH} 
\end{equation}
where $\mu$ is an unknown density, and
\begin{eqnarray}
  S \mu(\x) &=& \int_{\partial \Omega} G(\x,\y) \mu(\y) \, d\y \; ,\\ 
  D \mu(\x) &=& \int_{\partial \Omega} \partial_{n_y} G(\x,\y) \mu(\y) \, d\y
\end{eqnarray}
with $\partial_{n_y}$ representing the partial derivative in the direction normal to the boundary. $S \mu(\x)$ is known as a single layer potential, and $D \mu(\x)$ is known as a double layer potential \cite{Guenther88}.

For the harmonic function extension, we write $w$ as
\begin{equation}
  w(\x) = D \sigma (\x) + W \sigma
  \label{eq:layer_for_w}
\end{equation}
where $\sigma$ is an unknown density, and 
\begin{equation}
  W \sigma = \int_{\partial \Omega} \sigma(\y) \, d\y .
\end{equation}

Now, let $\mathcal{S} \mu (\x_0)$, $\mathcal{D} \mu (\x_0)$ and $\mathcal{D} \sigma (\x_0)$ denote the restrictions of $S$ and $D$ to points $\x_0$ on the 
boundary $\partial\Omega$, where the integrals are taken
in the Cauchy principal value sense when necessary. $u_H(\x)$ and $w(\x)$ 
reach the following limiting values as $\x$ approaches a point $\x_0$ on the boundary \cite{Guenther88}
\begin{equation}
\lim_{\x \to \x_0\, , \x \in \Omega} u^H(\x) = g(\x_0) - \tilde{v}(\x_0)=-\dfrac12 \mu(\x_0) + \mathcal{S} \mu(\x_0)
  + \mathcal{D} \mu(\x_0)
\label{eq:integral_mu}
\end{equation}
and
\begin{equation}
\lim_{\x \to \x_0\, , \x \in \Omega} w(\x) = f(\x_0)=\dfrac12 \sigma(\x_0)  + \mathcal{D} \sigma(\x_0) \; + W \sigma.
\label{eq:integral_sigma}
\end{equation}
\eqref{eq:integral_mu} is a second kind integral equation (SKIE) for $\mu$, and \eqref{eq:integral_sigma} an SKIE for $\sigma$. 

At this point, we have the desired integral representations for 
$u^H$ and for $w$, and equations for their associated densities. 
The representation \eqref{eq:layer_for_uH} for $u^H$ has 
been used in commercial software
\cite{MadMax02} and is known in the integral-equations
community \cite{Helsing05} but the authors are unaware of 
any treatment of the Fredholm alternative as applied to
the resulting integral equation \eqref{eq:integral_mu}.
We consider a proof of the invertibility of \eqref{eq:integral_mu} to be
beyond the scope of this paper but note that the argument of 
Lemma 29 in \cite{Rachh} can be modified to provide a proof
of its invertibility, even on multiply-connected domains.
The invertibility of \eqref{eq:integral_sigma} is well-known
\cite{Guenther88}.

We now discuss the numerical methods we chose to solve \eqref{eq:integral_mu} and \eqref{eq:integral_sigma} and to evaluate the integrals in \eqref{eq:layer_for_uH} and \eqref{eq:layer_for_w}.

\subsection{Solving the second kind integral equations for the densities}
In our solver, we discretize $\partial \Omega$ using panels of scaled, 16th
order Legendre nodes. Our numerical methods rely on the following simplifying assumptions concerning the boundary $\partial\Omega$: (1) the boundary is $C^k$ for some large $k$ and (2) the panels are chosen fine enough so that for source and target nodes on distinct, non-adjacent panels the integrals of the layer potentials are computed to high precision using the standard Gaussian weights (the ``source'' and ``target'' terminology is explained below). 
Note that there exist more complex algorithms that would allow us to relax both assumptions, and their implementation in our solver will be the subject of future work. To relax the first assumption, one could use 
any of the methods described in \cite{Helsing08C,Atkinson,Bremer10U,Serkh16}
to allow domains with corners.  While the second assumption is not necessarily much of a limitation on the types of domains which can be handled by our solver, the fineness implied by this assumption can lead to too great a computational burden for certain domains, such as domains in  which the boundary comes close to intersecting itself. The method of \cite{LocalQBX} provides a more efficient approach for such cases. 

Now, let $\partial \Omega$ be discretized into $L$ panels using  $M = 16L$ total nodes and denote the $i$th node by $\x_i$. Using generalized Gaussian quadrature for the interactions 
between nodes on the same and adjacent panels and the standard, scaled Gaussian weights otherwise, we obtain a Nystr\"{o}m discretization of \eqref{eq:integral_mu} and \eqref{eq:integral_sigma}:

\begin{eqnarray}
  g(\x_i)-\tilde{v}(\x_i) &=& -\dfrac12 \mu_i + 
   \sum_{j=1}^M \left ( G(\x_i,\x_j) \mu_j \omega^{s}_{i,j} +
  \partial_{n_j} G(\x_i,\x_j) \mu_j \omega^{d}_{i,j}\right) \; ,\label{eq:bdryintdisc}\\
  f(\x_i) &=& \dfrac12 \sigma_i + 
   \sum_{j=1}^M \left (
    \partial_{n_j} G(\x_i,\x_j) \sigma_j \omega^{d}_{i,j} +
   \sigma_j \omega_j \right )\label{eq:bdryintfunext}
\end{eqnarray}
where $\mu_i = \mu(\x_i)$, $\sigma_i = \sigma(\x_i)$ and the $\omega^s_{i,j}$, $\omega^d_{i,j}$, and $\omega_j$ correspond to integration weights. We note that the expressions above are a slight abuse of notation as the Green's function and its derivatives are undefined when $j = i$. The true formula is more generally a function of the boundary and the kernel but we find the above more edifying. In the current context, the relevant piece of information is that there exist weights $\omega^s_{i,j}$, $\omega^d_{i,j}$, and $\omega_j$ such that the quadratures appearing in the second kind equations can be evaluated with high-order accuracy. For a more detailed treatment of the generalized Gaussian quadrature framework, see \cite{Bremer10U}. In the following sections, we will refer to
the $\omega_j$, which are given by appropriately scaling the
standard Gauss-Legendre weights, as the smooth quadrature weights.

There exist many tools available for the fast solution of the linear systems \eqref{eq:bdryintdisc} and \eqref{eq:bdryintfunext}. There are iterative solution techniques, e.g. GMRES \cite{Saad},  
%which perform well when the eigenvalues of the system are clustered as they are for linear systems discretized from SKIEs
%and the system is well-conditioned. 
which perform well for linear systems discretized from SKIEs on
simple domains. The computational cost of such a scheme is  typically dominated by a term  of the order $kT$ where $k$ is the number of iterations required to converge and $T$ is the amount of work for a matrix-vector multiply.
For well-conditioned problems with $M$ boundary nodes, typically $k = \bigo (1)$ and the cost of $T$ can be reduced to $\bigo (M)$ with an FMM. There are also fast-direct solution methods, i.e., methods which construct,  in $\bigo(M)$ or $\bigo (M\log M)$ time, a
representation of the inverse of the system matrix which can be applied in $\bigo (M)$ or $\bigo (M\log M)$ time. For such direct methods, the cost of forming the representation
of the inverse is often much greater than that of the FMM, while the speed of applying the inverse, once computed, is often faster than the FMM. Fast-direct solvers can
be particularly useful for problems in highly-irregular domains, in which the iteration count of an iterative solver may be too high or unpredictable. They are also advantageous for cases in which several Laplace problems need to be solved for a fixed domain, since the high initial cost only has to be paid once. For our solver, we implemented the direct method developed by \cite{Ho12}, which is optimized for the type of problems considered here, and found that it gave very satisfactory performance. 

\subsection{Evaluation of $u^H$ and $w$ by quadrature-by-expansion}

Once $\mu$ and $\sigma$ are computed, we evaluate $u^H$ and $w$ by direct computation of the integrals \ref{eq:layer_for_uH} and \ref{eq:layer_for_w}. This step can at first appear complicated because the integral kernels are near singular for the evaluation of points near the boundary of the domain. We resolve the difficulty by computing the integrals for points near the boundary using the quadrature-by-expansion (QBX) method. We will not present the fundamentals of the QBX scheme here, since clear presentations for situations very closely related to the one we encounter here can be found in \cite{Kloeckner13,Ricketson,Askham2}. We will however stress two modifications to the standard QBX scheme which we implemented in our solver. First, we accelerated the evaluation of the layer potentials with the FMM (a similar but more sophisticated acceleration scheme is presented in \cite{Rachh16}). Second, we developed a variant of QBX which allows, after precomputation of the field at a fixed number of points, the evaluation of the field anywhere in the domain in $\bigo(1)$ time \cite{Askham2}. This is particularly convenient for the evaluation of the function extension when we construct the adaptive tree, since the grid points at which the values of the layer potential are desired may not be known a priori.

Let us discuss these two modifications to the 
standard QBX method in more detail. Consider, for 
example, the evaluation of the layer potential $u^H = S\mu + D\mu$. 
We use the notation as above for the discretization nodes 
$\mathbf{x}_{i}$, the smooth quadrature weights $\omega_i$, 
and the boundary normals $\mathbf{n}_i$ 
of $\partial\Omega$. Let $\mathbf{c}_{i}$ be the QBX centers,
located at a distance $r_i$ from the boundary: 
$\mathbf{c}_{i}=\mathbf{x}_{i}-r_i\mathbf{n}_i$. 
In the QBX method, the potential $u$ is approximated by a power 
series in the disc of radius $r_i$ about $\mathbf{c}_i$, denoted by
$B_{r_i}(\mathbf{c}_{i})$. For any $\mathbf{x}$ in $B_{r_i}(\mathbf{c}_{i})$, 
we write
\begin{equation}
u^H(\mathbf{x}) \approx \mbox{Re}\left(\sum_{l=0}^{p}\alpha_{l,i}(z-\xi)^{l}\right)
\; ,
\label{eq:power_series}
\end{equation}
where $z=x_{1}+ix_{2}$ and $\xi=c_{i,1}+ic_{i,2}$. The coefficients $\alpha_{l,i}$ can be recovered from the following integral on the circle of radius $r_i/2$ about $\mathbf{c}_{i}$:
\begin{equation}
\alpha_{l,i}=\frac{2^l}{\pi r_i^l}\int_{0}^{2\pi}u^H\left(\mathbf{c}_{i}+\frac{r_i}{2}(\cos\theta,\sin\theta)\right)e^{-il\theta} \; d\theta \; .
\label{eq:power_coeff}
\end{equation}

The layer potential $u^H$ is smooth on the circle of radius $r_i/2$, so the $\alpha_{l,i}$ can be computed with high order accuracy using the trapezoidal rule to discretize the integral \eqref{eq:power_coeff}. Let $M_{QBX}$ equispaced points $\mathbf{y}_{i,j}$ be placed on the circle $\partial B_{r/2}(\mathbf{c}_{i})$. The values $u^H(\mathbf{y}_{i,j})$ can be computed accurately using the smooth quadrature weights for $\partial \Omega$ to approximate the single and double layer potentials there, assuming that $r_i$ is large enough. For sufficient sampling, $M_{QBX}$ should be taken larger than $2p$. 

Once the coefficients are computed, the power series \eqref{eq:power_series} 
can be used to approximate the potential at targets which are close to
the boundary. High accuracy can be obtained when $r_i$ is sufficiently small.

With these preliminaries in place, the FMM-accelerated algorithm 
for the evaluation of the potential at $N$ targets $\mathbf{t}_{i}$ 
can be described in the following steps:

\begin{itemize}
\item Place $M$ centers at the points $\mathbf{c}_i = \mathbf{x}_i - r_i\mathbf{n}_i$.

\item Define $M_{QBX}$ equispaced points $\mathbf{y}_{i,j}$ for $j = 1,\ldots, M_{QBX}$ on the circle of radius $r_i/2$ about each center $\mathbf{c}_i$.

%\item Call the FMM with sources given by charges of strength $\mu(\mathbf{x}_i)\omega_i$ and dipoles oriented in the direction $\mathbf{n}_i$ and of the same strength located at the boundary points and targets given by the union of the $\mathbf{t}_i$ and the $\mathbf{y}_{i,j}$. This is a $\bigo(M M_{QBX} + N)$ procedure.

\item Call the FMM to evaluate $u^H$ at the targets $\mathbf{t}_i$ and
  the points $\mathbf{y}_{i,j}$, where the layer potentials are approximated using
  the smooth quadrature weights $\omega_i$. This is a $\bigo (MM_{QBX} + N)$ 
  procedure. 

\item Compute the coefficients $\alpha_{l,i}$ for each center as in \eqref{eq:power_coeff}, using the trapezoidal rule. This takes $\bigo (M M_{QBX} \log M_{QBX} )$ work for $M_{QBX} > p$ using the FFT.

\item For each target which is within $r_i/2$ of any boundary node $\x_i$, let $\mathbf{c}_j$ be the nearest QBX center. The smooth rule might not be accurate for this target, so
instead use the value given by the power series \eqref{eq:power_series} about $\mathbf{c}_i$. The cost for this is $\bigo(p)$ per target.
\end{itemize}

The scheme presented above is satisfactory if the targets $\mathbf{t}_{i}$ 
are all known in advance. However, when constructing the adaptive tree, 
the potential associated with the function extension may have to be evaluated 
at new targets $\mathbf{t}$. For the new targets which are close to the boundary, 
the potential can be computed using the expansion about the nearest QBX center. For the targets further from the boundary, we avoid calling the FMM again to compute the potential there. Instead, we store the multipole and local expansion coefficients computed for all boxes in the hierarchy during the initial call to the FMM. The values of the potential at the new targets $\mathbf{t}$ can then be evaluated in $\bigo(p_{FMM})$ work for each target, where $p_{FMM}$ is the order of the multipole and local expansions in the FMM \cite{Askham2}.

\begin{remark}
  In the above, we have avoided the key topic of how to select $r_i$, the radius of the QBX expansion. This must be done to balance two competing concerns: (1) that $r_i$ is sufficiently large so that $u^H(\mathbf{y}_{i,j})$ can be computed accurately using the smooth quadrature weights, $\omega_i$, and (2) that $r_i$ is sufficiently small so that the truncated power series \eqref{eq:power_series} is an accurate approximation of $u^H$. While it may seem unclear whether choosing an appropriate $r_{i}$ is indeed possible, this fact was proven in \cite{Epstein13}. Now, if the discretization node $\x_i$ is on a boundary panel of length $h_i$, it can be shown that setting $r_i= 4h_i$ provides high accuracy when using 16th order Legendre nodes on each panel. However, having the center so far places further restrictions on the discretization of the domain, since no boundary points are allowed to be in the interior of the QBX disc. Thus, in practice one often takes $r_i = h_i$ and computes $u^H(\mathbf{y}_{i,j})$ using the smooth weights for an oversampled version of the boundary.
\end{remark}

\begin{remark}
  The idea of using equispaced points on a circle to
  form a power series expansion of a harmonic function
  is reminiscent of the ``fast multipole method without
  multipoles'' of \cite{Anderson92}.
\end{remark}

\subsection{Derivatives of $u^H$}

In the QBX setting, the derivatives of the layer
potentials can be obtained by differentiating the 
QBX power series expansion. This works well for the single layer 
potential but has been observed to result in the loss
of precision for the double layer potential \cite{Kloeckner13}.
The cause of this loss of precision is unclear but may
result from the hyper-singular nature of the derivatives
of the double layer potential as operators on the 
boundary. 

It was pointed out to us by Manas Rachh and Leslie Greengard
\cite{IntByParts16} 
that the evaluation of the derivatives of the 
double layer potential can be accomplished by differentiating
the density along the boundary and using the QBX algorithm
for the Cauchy kernel (which behaves like the double
layer kernel). The key observations are that 
\begin{eqnarray}
  D\mu(\x) &=& \Re \left ( -\dfrac1{2\pi i} \int \dfrac{\mu(\xi)}{\xi - z} \, d\xi 
  \right ) \; , \\
  \left ( \partial_{x_1} - i \partial_{x_2} \right ) D\mu(\x) &=& 
  -\dfrac1{2\pi i} \int \dfrac{\mu(\xi)}{(\xi - z)^2} \, d\xi \; , \\
-\dfrac1{2\pi i} \int \dfrac{\mu(\xi)}{(\xi - z)^2} \, d\xi  &=& 
-\dfrac1{2\pi i} \int \dfrac{\mu'(\xi)}{\xi - z} \, d\xi \; ,
\end{eqnarray}
where $z = x_1 + i x_2$ and $\xi = y_1 + i y_2$ for $\y \in \Gamma$.
If high accuracy values of the tangential derivative of the
density can be obtained, the above provides a stable computation
of the derivatives of the double layer potential.
We have implemented this method but computing the tangential
derivative of the density presents some difficulty in the 
embedded boundary setting, as explained below.

In most of the numerical tests we perform in Section \ref{sec:results}, the discretization error, solution error, and QBX error associated with the boundary integral are made to be so small that the error in the boundary correction $u^H$ is dominated by the error in the interpolated values of the volume potential on the boundary. The accuracy of the gradient of the potential will be affected more strongly because the gradient of the  solution of the Laplace Dirichlet problem depends on the accuracy of the tangential derivative of the  boundary values. This is typically not a concern in the  integral equations context, for two reasons. First, one usually assumes that one  has high order accurate values for the boundary data. Second, the evaluation of the gradient of a layer
potential is smoothing for points sufficiently far from  the boundary. However, in the context of the Poisson solver we present in this article, the  tangential derivative of $\tilde{v}$ is, as an analytic matter, one order lower than the order of $\tilde{v}$ (this is because $\tilde{v}$ is a polynomial on each box) and the value  of the gradient may be requested arbitrarily close to (or even on) the boundary. This loss of accuracy for the boundary data results in a similar loss of accuracy for the computed gradient.

A potential way to address this problem is as follows. 
For a smoothly extended function $f_e$, a box code can  
be used to compute $\nabla V \tilde{f}_e$ at the collocation
points in the same way that $V \tilde{f}_e$ is computed.  
The interpolant of the gradient computed this way, which 
we call $\tilde{\mathbf{g}}$, is the same order as $\tilde{v}$.
We can then construct a new approximation to $V\tilde{f}_e$ on the boundary by first computing $\tilde{\mathbf{g}} \cdot \boldsymbol{\tau}$ along the boundary, where $\boldsymbol{\tau}$ is the tangent vector of the curve, and then computing its indefinite
integral panel-wise (we correct for the constant using $\tilde{v}$, again panel-wise). The resulting function is the same order accuracy as $\tilde{v}$ but its derivative
is a better approximation to the derivative of $V\tilde{f}_e$. We investigate the merits of this alternative approach numerically in Section \ref{sec:results}.

Unfortunately, this approach does not appear to improve the order
of accuracy for non-smooth $f_e$. As noted in Section \ref{subsec:nonsmoothf},
we expect the convergence of the computed gradient to be one 
order lower than the potential for an extended density $f_e$
obtained through continuous extension or extension by zero, i.e.
there is no advantage to using $\tilde{\mathbf{g}}$ as it is already one
order of accuracy lower than $\tilde{v}$. 
This was evident in numerical experiments,
where the accuracy in the derivatives was comparable using
this new approximation to $V\tilde{f_e}$ instead of $\tilde{v}$.
We therefore leave these results out of the next section.
In the case that a $C^k$ extension is available for sufficiently 
large $k$, this technique may prove important to achieving 
super-convergence.

\section{Numerical results}
\label{sec:results}

In order to verify the preceding analysis and  test the performance of the numerical method we propose in this article, we have implemented a Poisson solver in Fortran which combines all the different modules we presented in the previous sections. The volume integral code is a modified version of the original Fortran code of \cite{Ethridge} (using fourth order nodes on leaf boxes), with some added OpenMP parallelism and the modification for computing gradient values discussed in Section \ref{subsec:volderivatives}. The codes for the boundary correction and the continuous extension were written specifically for this work, and are based on the methods described in Section \ref{sec:irregular_box} and \ref{sec:layer}. We are currently documenting the numerical solver we used to generate the results shown below, and will make it freely available online at a later date.

\begin{figure}[h!] 

\centering

\includegraphics[width=.8\textwidth]{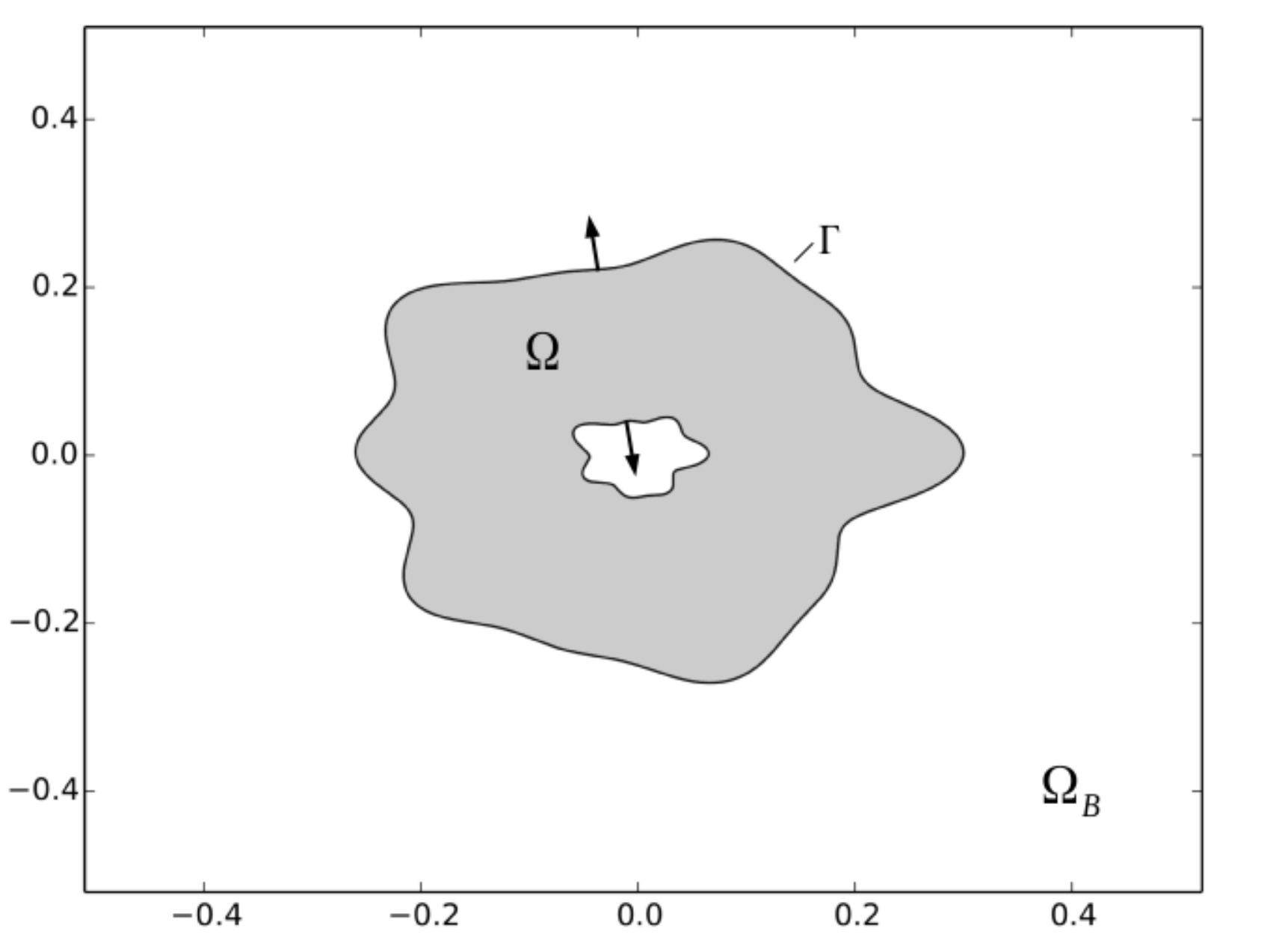}

\caption{The domain $\Omega$ and its boundary $\Gamma$. The
axes coincide with the boundary of the containing box $\Omega_B$.
Two outward-pointing boundary normal vectors are indicated by
arrows.}

\label{fig:omega}

\end{figure}

For each numerical test, we have used the domain $\Omega$ shown in Figure \ref{fig:omega},
which has an irregular boundary and is multiply connected. The interfaces of this
domain are specified by parametric equations in polar coordinates. Specifically, 
each interface is given by a set of points $(\theta,r(\theta))$ for $\theta \in [0,2\pi)$,
where $r(\theta) = c_0 + \sum_j (c_j \cos(j\theta) + d_j\sin(j\theta))$. The choice
of the coefficients was arbitrary. For reference, the non-zero coefficients for the
outer boundary were $c_0 = 0.25$, $d_3 = c_6 = c_8 = c_{10} = 0.01$, $c_5 = 0.02$. 
The non-zero coefficients for the inner boundary were $c_0 = 0.05$, 
$c_2 = d_3 = c_5 = c_7 = 0.005$.

Let $N_p$ denote the number of panels used in the discretization of the boundary and $M=16*N_p$ denote the total number of boundary points (we use 16th order nodes throughout). For the  volume integral nodes, let $N_V$ denote the total number of points in $\Omega_B$ and 
$N_\Omega$ denote the number of points inside $\Omega$. For each test, we approximated the relative $L_\infty$ error,
$$ E(\psi) = \dfrac{\displaystyle\max_\Omega |\psi_{exact}-\psi_{numerical}|}{\displaystyle\max_\Omega |\psi_{exact}|}\; ,$$ 
where $\psi$ is either the potential or its derivatives, by sampling at $10^6$ 
points randomly placed in $\Omega$. These points were kept the same for each discretization level for the sake of convergence tests. We report the error in the gradient 
below as $\sqrt{E(u_x)^2 + E(u_y)^2}$. In the error analysis of this section, the density for the boundary correction is computed with high accuracy (say 12 digits) and the corresponding layer potential is evaluated with high accuracy as well (say 
12 digits for the potential and 8-9 digits for its gradient). With this assumption, the error will be  primarily a function of the number of discretization nodes in the volume, i.e. $N_\Omega$.

All computations were performed on a desktop computer with an Intel® Xeon(R) CPU E3-1220 v5 (3.00GHz, 4 core) and 16 Gb of memory. A few of the computations depend only on the boundary and therefore take the same amount of time for each discretization level. In the first example, the boundary was discretized with $M=9,280$ nodes. The precomputation time for the direct solver
took 1.20 and 1.80 seconds for the continuous function
extension and boundary correction linear systems, 
respectively. The precomputation time to allow for 
$\bigo (1)$ access to the layer potential took 
.60 and 1.20 seconds for the continuous function 
extension and boundary correction layer potentials, 
respectively. The solution in the second example is 
much more irregular than in the first and thus more 
boundary points were required. For this case, the 
boundary was discretized with $M=14,208$
nodes. The precomputation time for the direct solver
took 1.57 and 2.89 seconds for the continuous function
extension and boundary correction linear systems, 
respectively. The precomputation time to allow for 
$\bigo (1)$ access to the layer potential took 
1.10 and 1.99 seconds for the continuous function 
extension and boundary correction layer potentials, 
respectively. We consider these computational costs
to be modest and emphasize that the precomputation
of the fast-direct solver must only be done once
per domain. We report on the speed of the volume 
integral code and the evaluation of the layer potentials
below.

\subsection{A note on adaptivity}

When a smooth extension $f_e$ is known,
the bound \eqref{eq:errorAnalysis} of Section
\ref{subsec:smoothf} implies a rather straightforward
\textit{a priori} adaptive discretization strategy: 
for a given tolerance, refine the tree until the local polynomial
interpolant on each leaf box approximates $f_e$ 
within that tolerance, which can be tested by comparing
$f_e$ and the interpolant on a finer grid. It turns out that this strategy will result in an 
overall error well below the desired tolerance. A modification which, in practice, gets closer to the desired tolerance is to refine until the error in the local polynomial interpolation times the area of the box is within the tolerance on each leaf box. 

For piecewise smooth $f_e$, we saw in Section \ref{subsec:nonsmoothf} that the bound 
\eqref{eq:errorAnalysis} may be pessimistic. However, the analysis of that section offers
little in terms of an \textit{ a priori} discretization strategy. If the above strategy for smooth $f_e$ is implemented, the accuracy of the resulting scheme is often not even competitive with a uniform grid. We consider the problem of efficient \textit{a priori} adaptive discretization to be open in this setting but have empirically found the following scheme to compare favorably to uniform discretization in our tests: weight the error approximation using the size of the given leaf box as described above
but using the area of the leaf box for boxes 
which intersect the boundary (where $f_e$ is less 
smooth) and using the sidelength of the
leaf box otherwise. In this sense, we are more
forgiving of the approximation error for boxes
where $f_e$ is merely continuous. 

It seems that in many situations an \textit{a posteriori} discretization strategy would be more efficient in terms of accuracy per grid point. While this may be an intuitive statement, it is not clear whether an \textit{a posteriori} scheme would be more efficient in terms of 
accuracy per flop because such a scheme may require several successive iterations. We do not attempt to answer this question here but emphasize that the issue with \eqref{eq:errorAnalysis} is a matter of efficiency rather than correctness, i.e. if an \textit{a priori} bound is required, \eqref{eq:errorAnalysis} provides one, it just may be an over-estimate.

\subsection{Example 1}

\begin{figure}[h!]
  \begin{centering}
    \includegraphics[width=.9\textwidth]{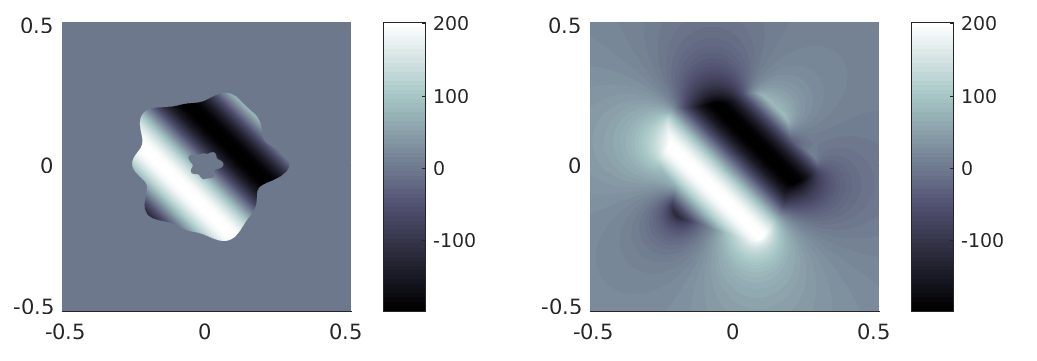}
    \caption{The extended density $f_e$ for Example 1 
      using extension by zero (left) and continuous extension (right).}
  \label{fig:ex1ext}
  \end{centering}
\end{figure}

For Example 1, we choose a known, relatively smooth solution $u$ given by  

\begin{equation}
  u(\mathbf{x}) = \sin(10(x_1+x_2)) + x_1^2 - 3x_2 + 8 \;.
  \label{eq:example1}
\end{equation}
and calculate $f$ analytically by direct differentiation, and $g$ by evaluating $u$ on $\partial\Omega$. Figure \ref{fig:ex1ext} shows heat maps of the corresponding $f_e$ obtained by zero extension and by continuous extension. Example 1 is relatively simple on purpose, in order to test the validity of the analysis of the previous sections.

First, consider the question of superconvergence for a smooth extension $f_e$. This is simple to test numerically as the formula \eqref{eq:example1} for $u$ is smooth on $\R^2$. In Section \ref{sec:layer}, we noted that the boundary correction can be computed with two different types of boundary data. Let version 1 denote the boundary data obtained from $\tilde{v}$ and version 2 denote the boundary data obtained by integrating $\boldsymbol{\tau} \cdot \tilde{\mathbf{g}}$, where we have reused the notation of Section \ref{sec:layer}. We perform a convergence test on uniform trees for both versions 1 and 2. 
According to the analysis of the preceding sections, version 1 should display fourth order convergence for the potential and sub-fourth order convergence for the gradient, while version 2 should display superconvergence, i.e. fourth order for both the potential and gradient. 

\begin{figure}[h!]
  \begin{centering}
    \includegraphics[width=.45\textwidth]{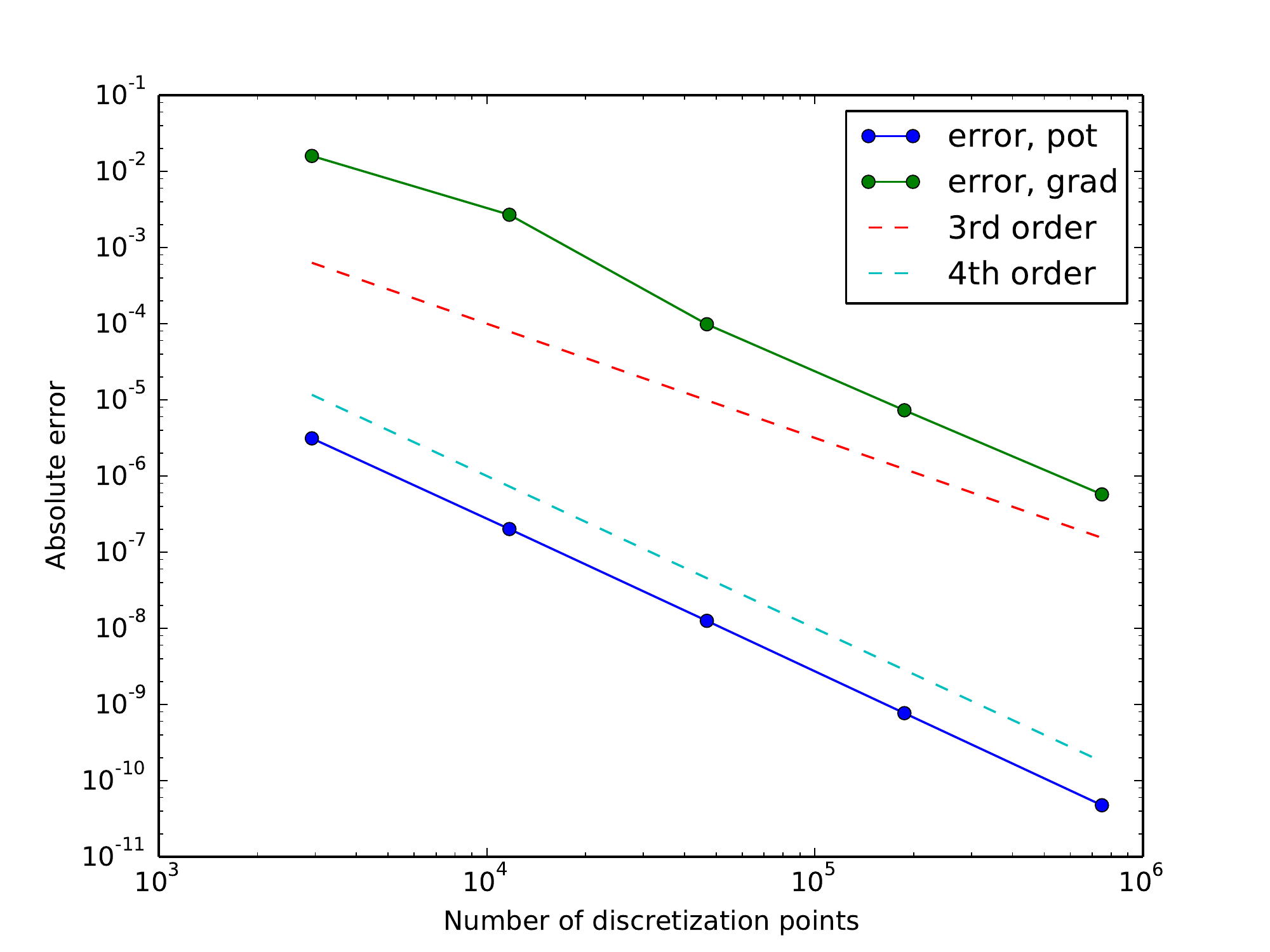}
    \includegraphics[width=.45\textwidth]{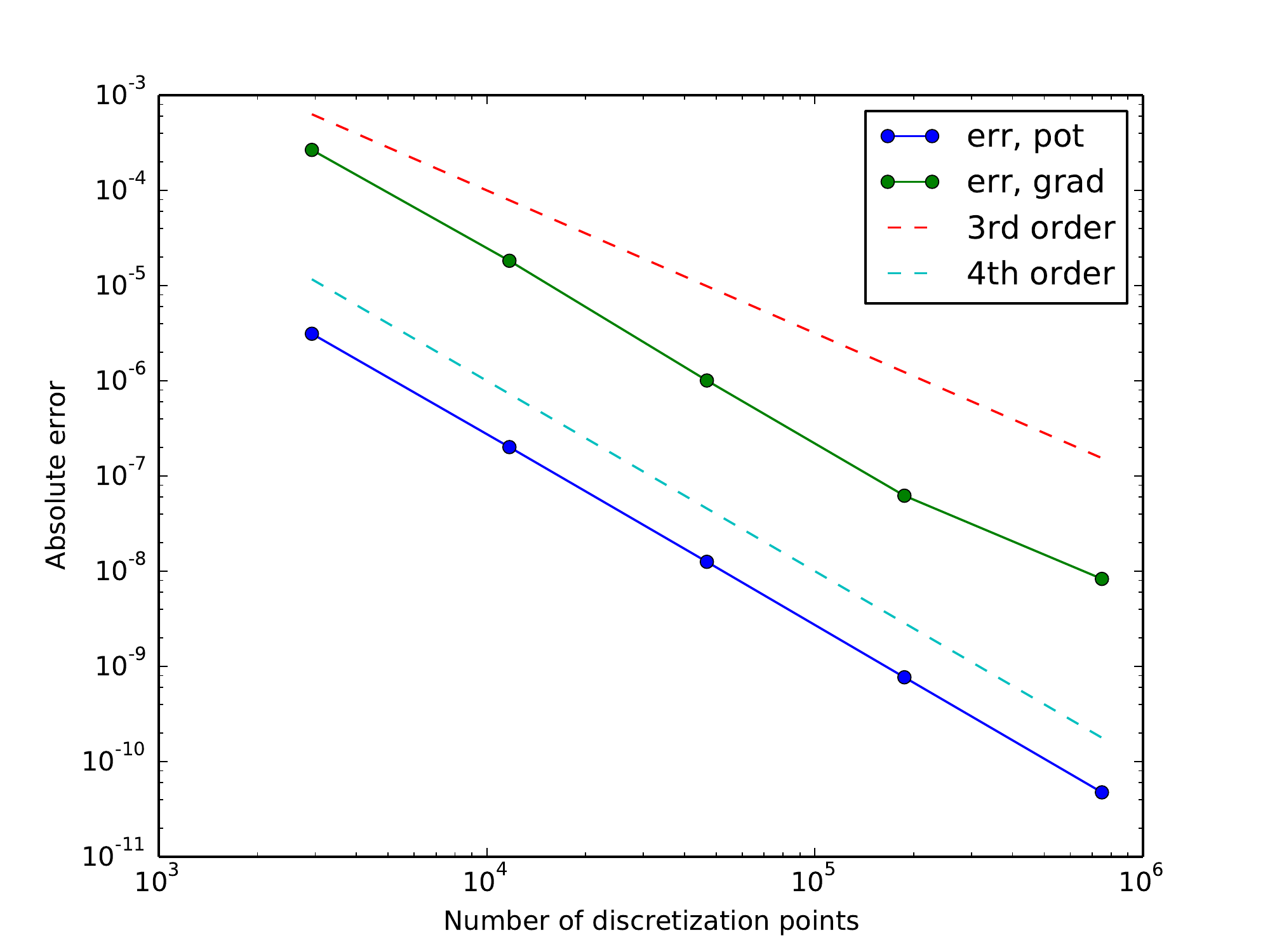}
    \caption{Example 1, smooth extension. 
      Accuracy of the potential and gradient versus the number of discretization nodes $N_\Omega$ for version 1 (left) and version 2 (right).}
  \label{fig:ex1smooth}
  \end{centering}
\end{figure}

In Figure \ref{fig:ex1smooth}, we see that the analysis is largely confirmed. While we cannot conclude decisively regarding the convergence order of the gradient for version 1, it is indeed fourth order for version 2. Note that the slope seems to taper off for the last point, which is likely due to the fact that one is approaching the accuracy of the QBX evaluation of the derivative. In terms of accuracy per grid point, version 2 is clearly superior to version 1.

Next, we consider the question of the convergence order using extension-by-zero and continuous extension with a  layer potential. The analysis of Section \ref{subsec:nonsmoothf}
suggests that we should see second order convergence for the potential and first order convergence for the gradient using extension-by-zero. This should be improved to third order for the potential and second order for the gradient by using continuous extension. As a reminder, these rates are to be compared with the rates implied by the coarser error bound \eqref{eq:errorAnalysis}, which suggests that the extension-by-zero scheme would not converge and that the continuous extension  scheme would be merely first order in the potential and derivative. To test the reasoning of Section \ref{subsec:nonsmoothf}, we performed a convergence test of the extension by zero and continuous extension methods on uniform trees.

\begin{figure}[h!]
  \begin{centering}
    \includegraphics[width=.45\textwidth]{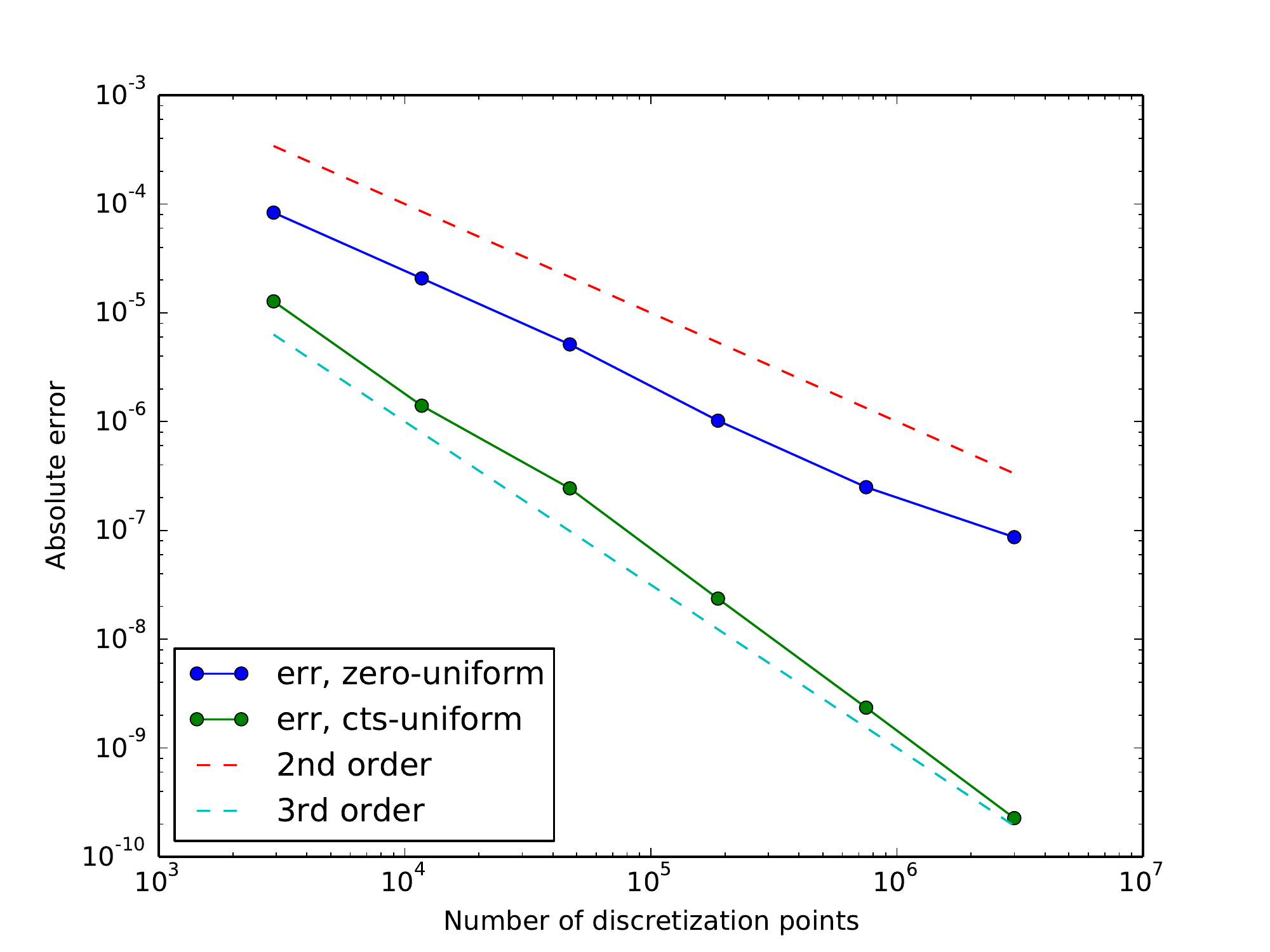}
    \includegraphics[width=.45\textwidth]{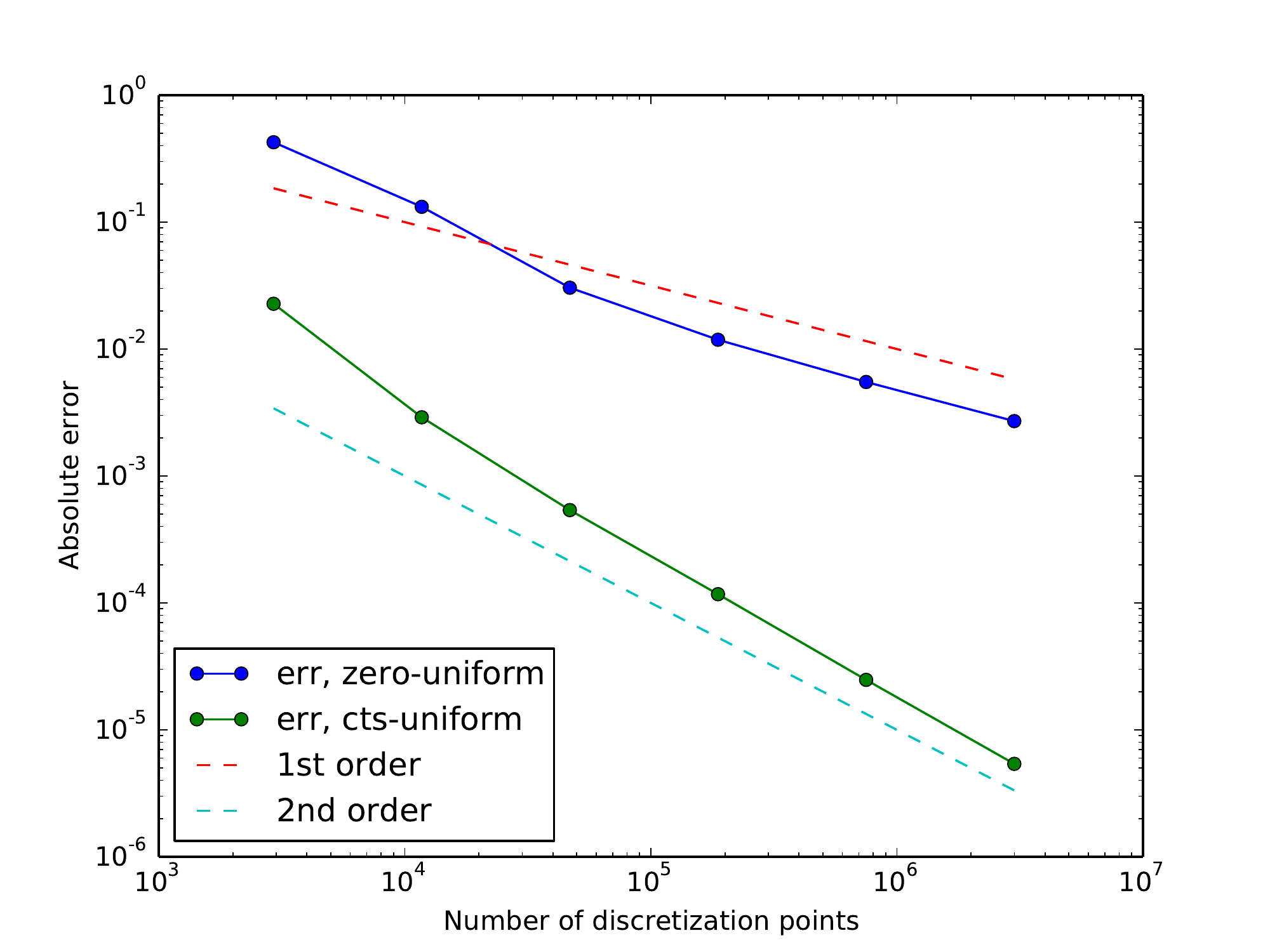}
    \caption{Example 1, convergence rates
      on a uniform tree. Accuracy of the potential (left) and gradient (right) versus the number of discretization nodes $N_\Omega$, using either extension-by-zero (blue curves) or continuous extension (green curves).}
  \label{fig:ex1uni}
  \end{centering}
\end{figure}

The results are shown in Figure \ref{fig:ex1uni}, and confirm that the  analysis of Section \ref{subsec:nonsmoothf} gives a better sense of the convergence rate than a na\"{i}ve application of the bound \eqref{eq:errorAnalysis}.

\begin{figure}[h!]
  \begin{centering}
    \includegraphics[width=.9\textwidth]{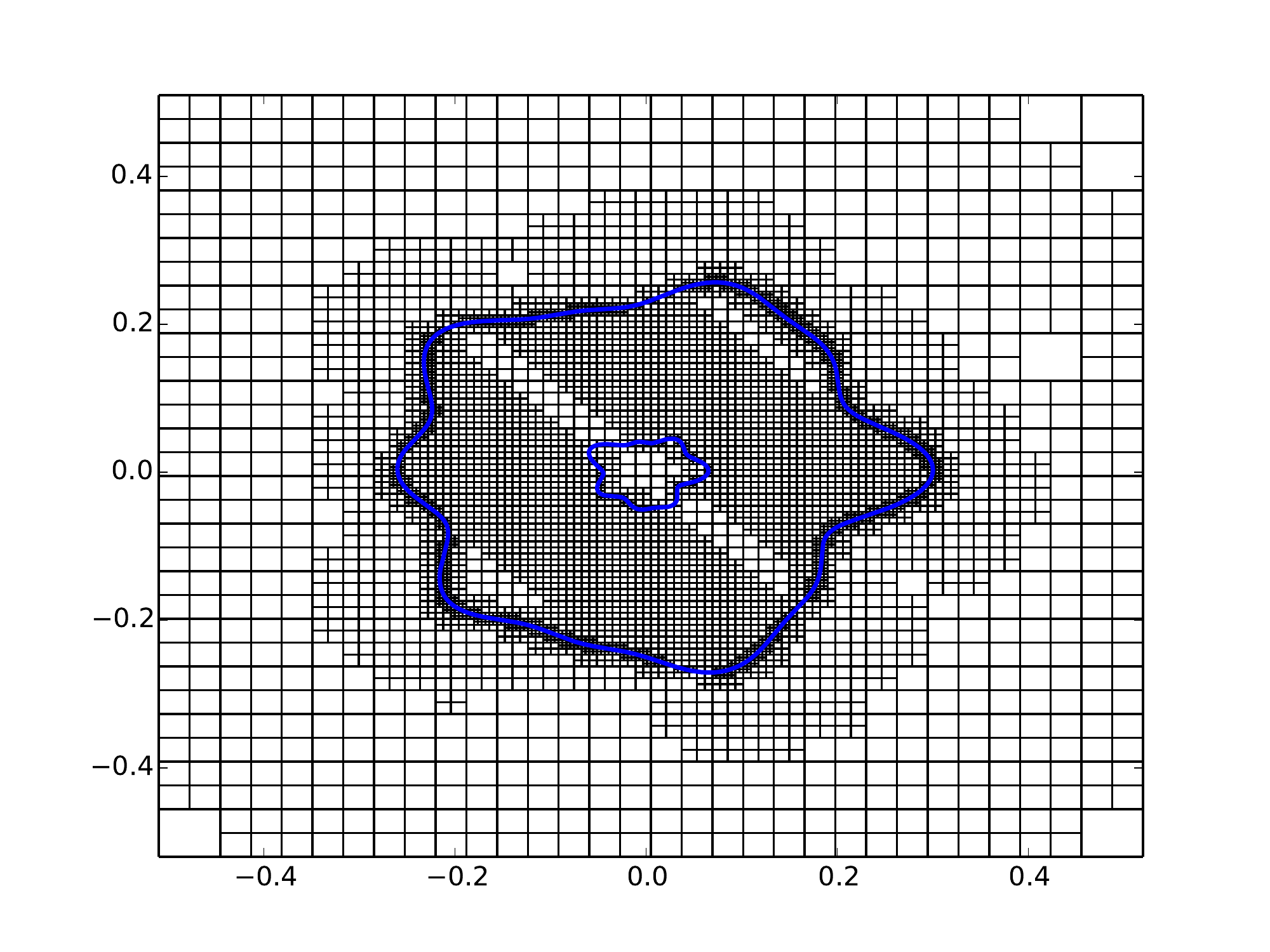}
    \caption{Example 1. An example adaptive tree for 
      the continuously extended $f_e$.}
  \label{fig:ex1adpboxes}
  \end{centering}
\end{figure}

Next, we consider the question of adaptive grid refinement. An adaptive grid should be able to provide significant gains, especially for the nonsmooth $f_e$. For the results presented here, we use an adaptive tree based on \textit{a priori} error estimates, as described 
in the previous subsection. This refinement rule tends to place more boxes near the boundary because  of the irregularity of $f_e$ across the boundary, as shown in Figure \ref{fig:ex1adpboxes}. We only present results corresponding to continuous extension here, as our refinement rule did not work well with zero extension, and is not relevant in the case of the smooth $f_e$ since there is little difference between adaptive and uniform discretization in that case.

\begin{figure}[h!]
  \begin{centering}
    \includegraphics[width=.45\textwidth]{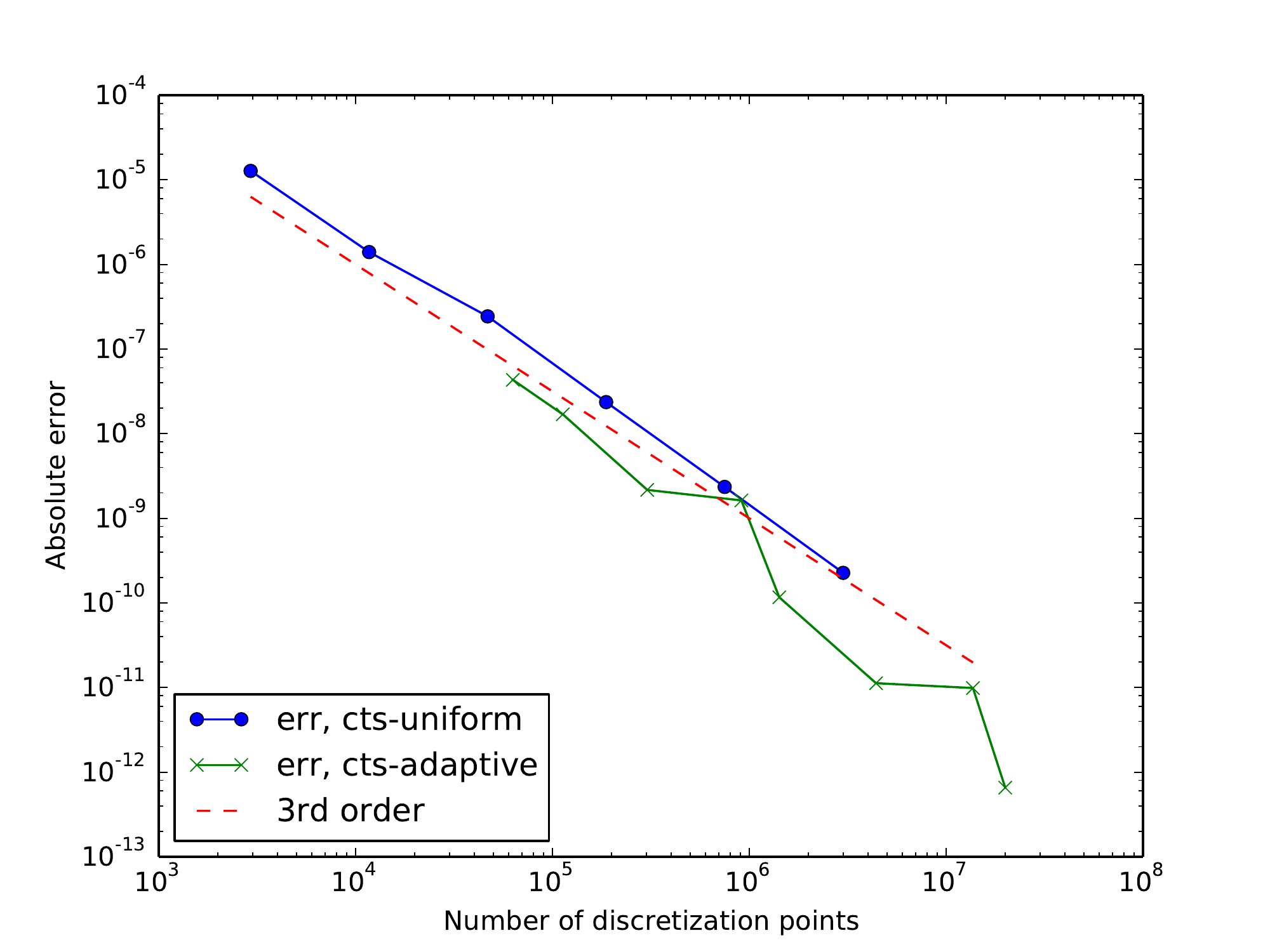}
    \includegraphics[width=.45\textwidth]{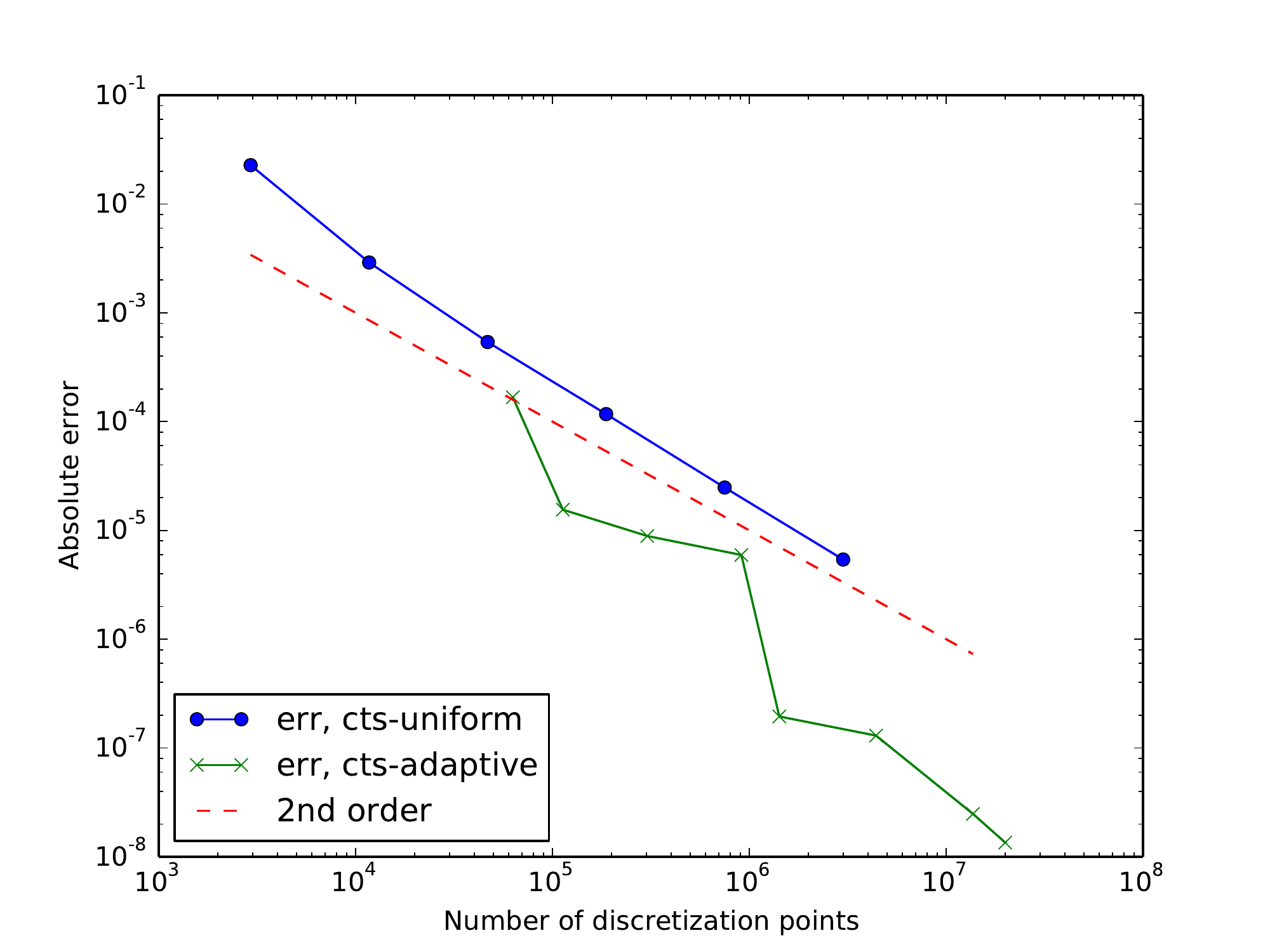}
    \caption{Example 1. Using continuous extension,  
      a plot of the error in the potential (left) and error 
      in the gradient (right) versus the number of discretization 
      nodes $N_\Omega$,
      for both uniform (blue dots) and adaptive trees (green crosses).}
  \label{fig:ex1adp}
  \end{centering}
\end{figure}

In Figure \ref{fig:ex1adp}, we see modest improvement in
the accuracy of the potential and significant improvement
in the accuracy of the gradient using adaptive discretization.
We note that for the tests with adaptive grids
much larger values of $N_\Omega$ could be achieved. This is because the memory 
consumption of the volume integral code depends on $N_V$, 
the total number of nodes in the box $\Omega_B$. The
uniform tree rather inefficiently places many points
outside of $\Omega$, whereas the adaptive tree places
relatively few points because the extended function is quite
smooth outside of $\Omega$, where it is harmonic.

We conclude this section on Example 1 by analyzing the run time performance of the box code and of the evaluation of layer potentials. Figures for the box code are given in Tables \ref{tab:ex1unitimes} and \ref{tab:ex1adptimes}, and figures for the evaluation of layer potentials are given in Tables \ref{tab:ex1uniqbxtimes} and \ref{tab:ex1adpqbxtimes}. $t_V$ denotes the time for the box code, $t_{QP}$ denotes the time for QBX precomputation (forming the expansions for $\bigo (1)$ access to the field, as described above), and $t_{QE}$ denotes the time for QBX evaluations at each node in the domain. Each of these times includes the time required to evaluate both the potential and the gradient. The performance is only reported for continuous extension; the results for extension by zero and smooth extension are similar.

\begin{table}
\centering
\begin{tabular}{lllll} 
\toprule 
$N_\Omega$ & $N_V$ & $t_V$ & $N_\Omega/t_V$ & $N_V/t_V$ \\ 
\midrule 
2.9290e+03 & 1.6384e+04 & 3.0193e-02 & 9.7009e+04 & 5.4264e+05 \\ 
1.1717e+04 & 6.5536e+04 & 4.4461e-02 & 2.6353e+05 & 1.4740e+06 \\ 
4.6846e+04 & 2.6214e+05 & 1.2574e-01 & 3.7256e+05 & 2.0848e+06 \\ 
1.8739e+05 & 1.0486e+06 & 3.9343e-01 & 4.7629e+05 & 2.6652e+06 \\ 
7.4955e+05 & 4.1943e+06 & 1.4926e+00 & 5.0218e+05 & 2.8101e+06 \\ 
2.9983e+06 & 1.6777e+07 & 7.3144e+00 & 4.0991e+05 & 2.2937e+06 \\ 
\bottomrule 
\end{tabular}
\caption{Box code timing information for Example 1 with continuous
function extension and a uniform tree.\label{tab:ex1unitimes}}
\end{table}

\begin{table}
\centering
\begin{tabular}{lllll} 
\toprule 
$N_\Omega$ & $N_V$ & $t_V$ & $N_\Omega/t_V$ & $N_V/t_V$ \\ 
\midrule 
6.2928e+04 & 1.1685e+05 & 7.5608e-02 & 8.3229e+05 & 1.5454e+06 \\ 
1.1291e+05 & 2.3666e+05 & 1.2941e-01 & 8.7251e+05 & 1.8287e+06 \\ 
3.0310e+05 & 5.6781e+05 & 2.7364e-01 & 1.1077e+06 & 2.0750e+06 \\ 
9.1144e+05 & 1.5124e+06 & 6.4442e-01 & 1.4144e+06 & 2.3468e+06 \\ 
1.4207e+06 & 2.7318e+06 & 1.1932e+00 & 1.1906e+06 & 2.2895e+06 \\ 
4.4043e+06 & 7.3749e+06 & 3.0303e+00 & 1.4534e+06 & 2.4337e+06 \\ 
\bottomrule 
\end{tabular}
\caption{Box code timing information for Example 1 with continuous
function extension and an adaptive tree.
\label{tab:ex1adptimes}}
\end{table}

\begin{table}
\centering
\begin{tabular}{lllll} 
\toprule 
$N_\Omega$ & $t_{QP}$ & $t_{QE}$ & $N_\Omega/(t_{QP}+t_{QE})$ & $N_\Omega/t_{QE}$ \\ 
\midrule 
2.9290e+03 & 1.1848e+00 & 1.3120e-03 & 2.4694e+03 & 2.2325e+06 \\ 
1.1717e+04 & 1.1697e+00 & 3.2675e-03 & 9.9892e+03 & 3.5859e+06 \\ 
4.6846e+04 & 1.1825e+00 & 1.2799e-02 & 3.9192e+04 & 3.6601e+06 \\ 
1.8739e+05 & 1.2034e+00 & 5.4077e-02 & 1.4902e+05 & 3.4652e+06 \\ 
7.4955e+05 & 1.1677e+00 & 1.8898e-01 & 5.5249e+05 & 3.9663e+06 \\ 
2.9983e+06 & 1.1896e+00 & 7.4644e-01 & 1.5487e+06 & 4.0168e+06 \\ 
\bottomrule 
\end{tabular}
\caption{QBX timing information for Example 1 with 
continuous function extension and a uniform tree.
\label{tab:ex1uniqbxtimes}}
\end{table}

\begin{table}
\centering
\begin{tabular}{lllll} 
\toprule 
$N_\Omega$ & $t_{QP}$ & $t_{QE}$ & $N_\Omega/(t_{QP}+t_{QE})$ & $N_\Omega/t_{QE}$ \\ 
\midrule 
6.2928e+04 & 1.1912e+00 & 1.7972e-02 & 5.2042e+04 & 3.5014e+06 \\ 
1.1291e+05 & 1.1744e+00 & 3.1084e-02 & 9.3665e+04 & 3.6325e+06 \\ 
3.0310e+05 & 1.1708e+00 & 7.9357e-02 & 2.4245e+05 & 3.8195e+06 \\ 
9.1144e+05 & 1.1733e+00 & 2.3131e-01 & 6.4889e+05 & 3.9404e+06 \\ 
1.4207e+06 & 1.1770e+00 & 3.6280e-01 & 9.2264e+05 & 3.9159e+06 \\ 
4.4043e+06 & 1.1733e+00 & 1.1130e+00 & 1.9264e+06 & 3.9572e+06 \\ 
\bottomrule 
\end{tabular}
\caption{QBX timing information for Example 1 with 
continuous function extension and an adaptive tree.
\label{tab:ex1adpqbxtimes}}
\end{table}

There are a few points to highlight from Tables \ref{tab:ex1unitimes} and \ref{tab:ex1adptimes}. We see that $N_V/t_V$ is
roughly constant for large $N_V$, indicating that the 
FMM indeed scales linearly in terms of the total 
number of FMM nodes. One of the strengths of a box code
is that this ratio is similar for uniform and adaptive 
tres. Further, the throughput is quite good,
at about 2.5 million points per second. We include
the ratio $N_\Omega/t_V$ because the number of grid points 
inside the domain seems to be the more natural figure of merit. 
For a uniform tree (Table \ref{tab:ex1unitimes}), we have
that $N_\Omega$ is a fixed fraction of $N_V$, so that 
$N_\Omega/t_V$ is some fraction of $N_V/T_V$; here it is
typically around 470 thousand points per second. In the
adaptive case (Table \ref{tab:ex1adptimes}), the nodes
can be placed more intelligently inside the domain
and we see that the throughput --- in terms of $N_\Omega/t_V$ ---
is better than in the uniform case. 

Tables \ref{tab:ex1uniqbxtimes} and \ref{tab:ex1adpqbxtimes},
show that the run time for QBX is similar for volume nodes arranged in uniform or adaptive trees, as one might expect. If one only considers the cost of the evaluations, we see that the throughput, $N_\Omega/t_{QE}$, is roughly constant at about 3.9 million points per second. The precomputation time, $t_{QP}$,  depends only on the number of boundary nodes $M$ and is large  relative to $t_{QE}$ until $N_\Omega$ is of the order of a few millions. 
When this precomputation time is included, the throughput,
$N_\Omega/(t_{QP}+t_{QE})$, is still quite high, on the same
order as the box code for large $N_\Omega$. Of course, for a boundary with many
discretization nodes $M$, one expects this to no longer be the case.

\subsection{Example 2} 

For Example 2, we choose an exact solution $u$ with a sharp ridge along the $x_2$ axis, given by
\begin{equation}
  u(\mathbf{x}) = \sin(10(x_1+x_2)) + x_1^2 - 3x_2 + 8 + e^{-500x_1^2} \;.
\end{equation}
As before, we obtain a closed form formula for $f$ by calculating the Laplacian of $u$. Observe that $f$ has very sharp variations. This example was chosen on purpose to specifically illustrate and analyze the value of adaptive mesh refinement. As in Example 1, $g$ is computed with arbitrary accuracy by evaluating $u$ on $\partial\Omega$. The function $g$ also has sharp variations, and so does the volume integral. In order to better resolve the boundary data we thus use $M=14,208$ boundary nodes in this example, as opposed to $M=9,280$ in Example 1. 

\begin{figure}[h!]
  \begin{centering}
    \includegraphics[width=.45\textwidth]{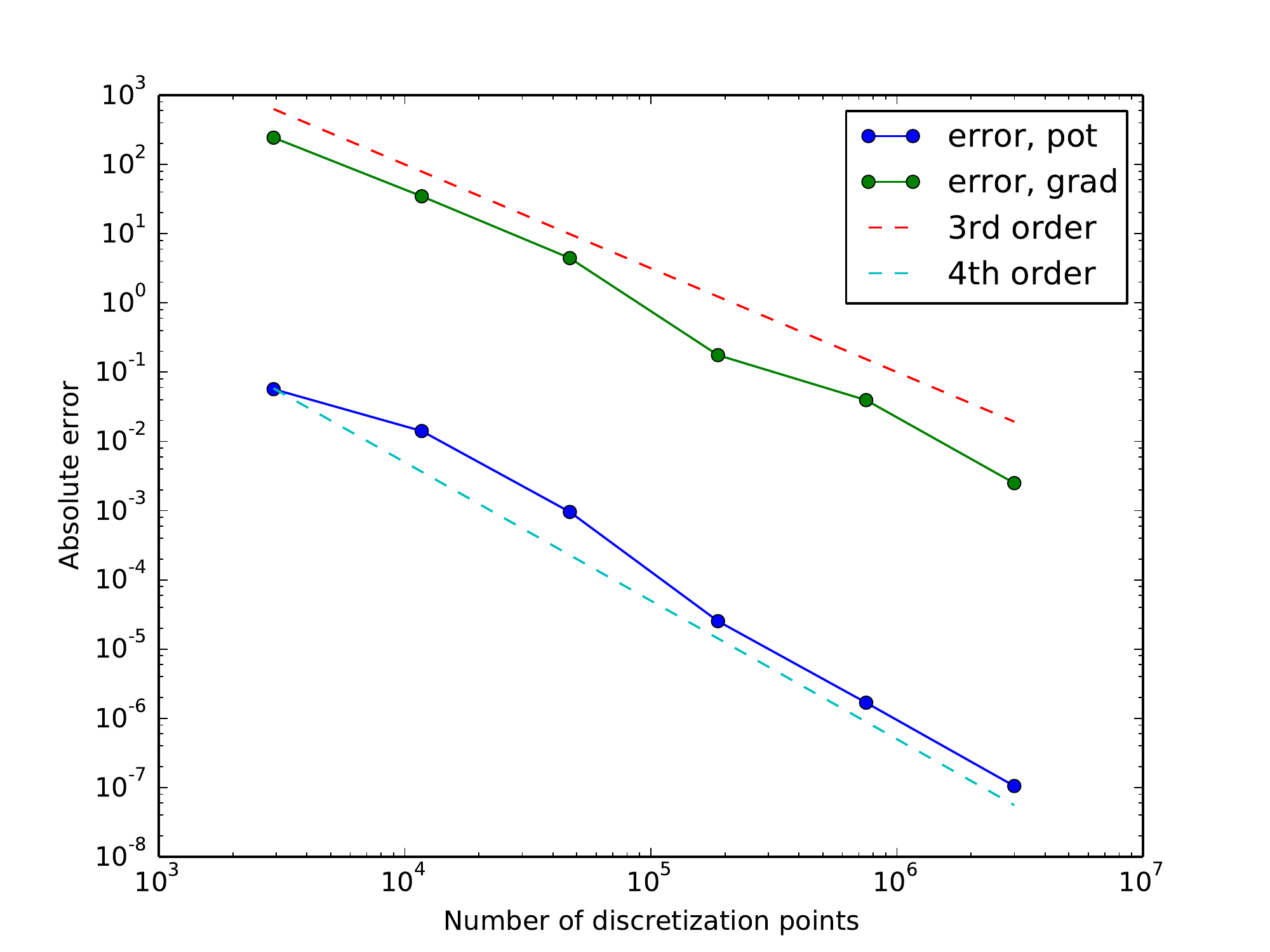}
    \includegraphics[width=.45\textwidth]{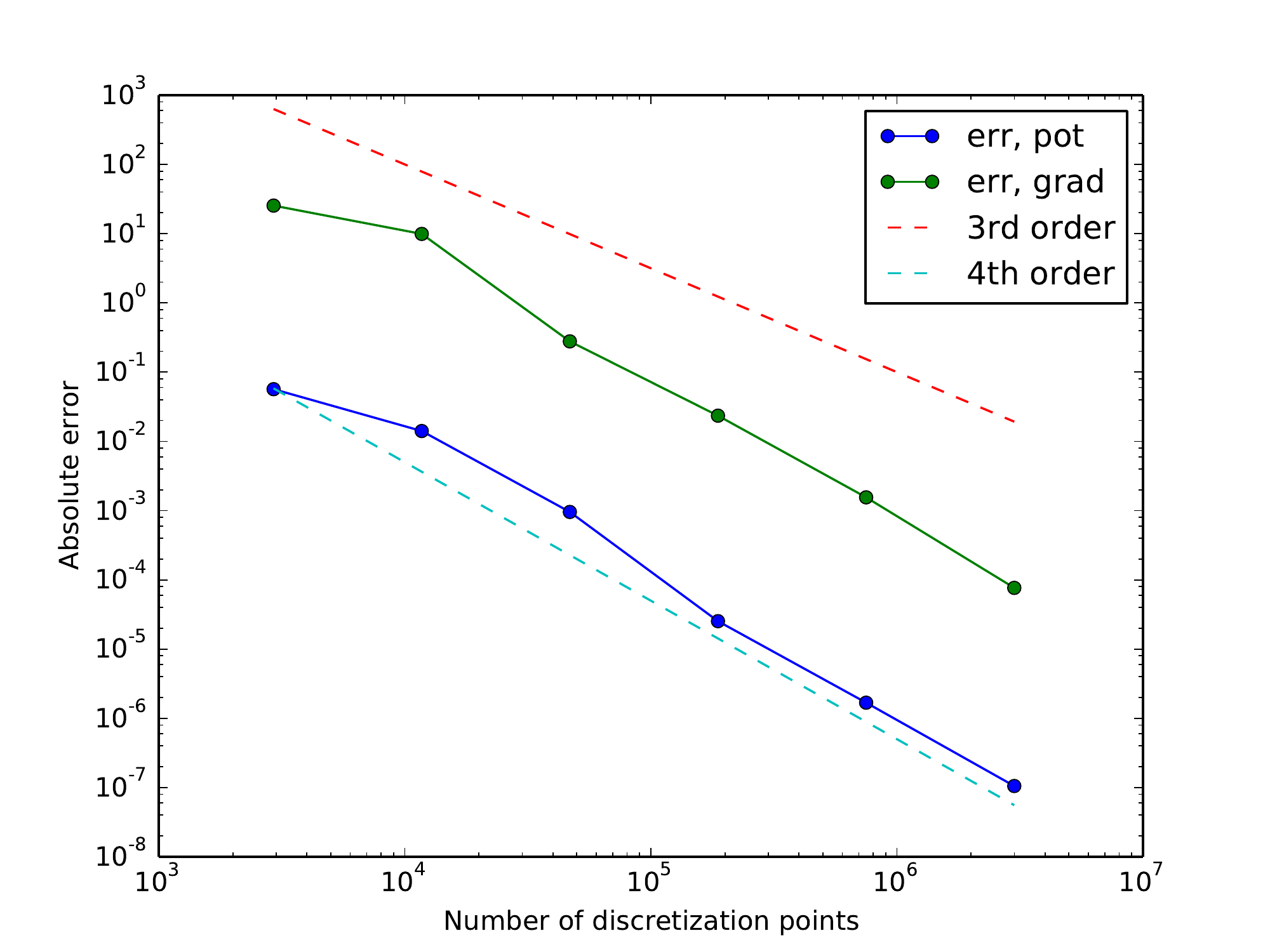}
    \caption{Example 2, smooth extension. 
      Accuracy of the potential and gradient versus the number of discretization nodes $N_\Omega$
      for version 1 (left) and version 2 (right).}
  \label{fig:ex2smooth}
  \end{centering}
\end{figure}

First, consider the question of superconvergence for a smooth extension $f_e$.
Let version 1 and version 2 of the boundary data be defined as in Example 1. We perform a convergence test on uniform trees for both versions 1 and 2. 
As before, version 1 should display fourth order convergence for the potential and sub-fourth order convergence for the gradient, while version 2 should display superconvergence. This is precisely what we see in Figure \ref{fig:ex2smooth}. For each version, the initial
convergence order is slow, likely a result of the irregularity of $f$. It is unclear what the  eventual convergence order of the gradient is for version 1 but it is fourth order for version 2. As in Example 1, the accuracy of version 2 is much better.

\begin{figure}[h!]
  \begin{centering}
    \includegraphics[width=.45\textwidth]{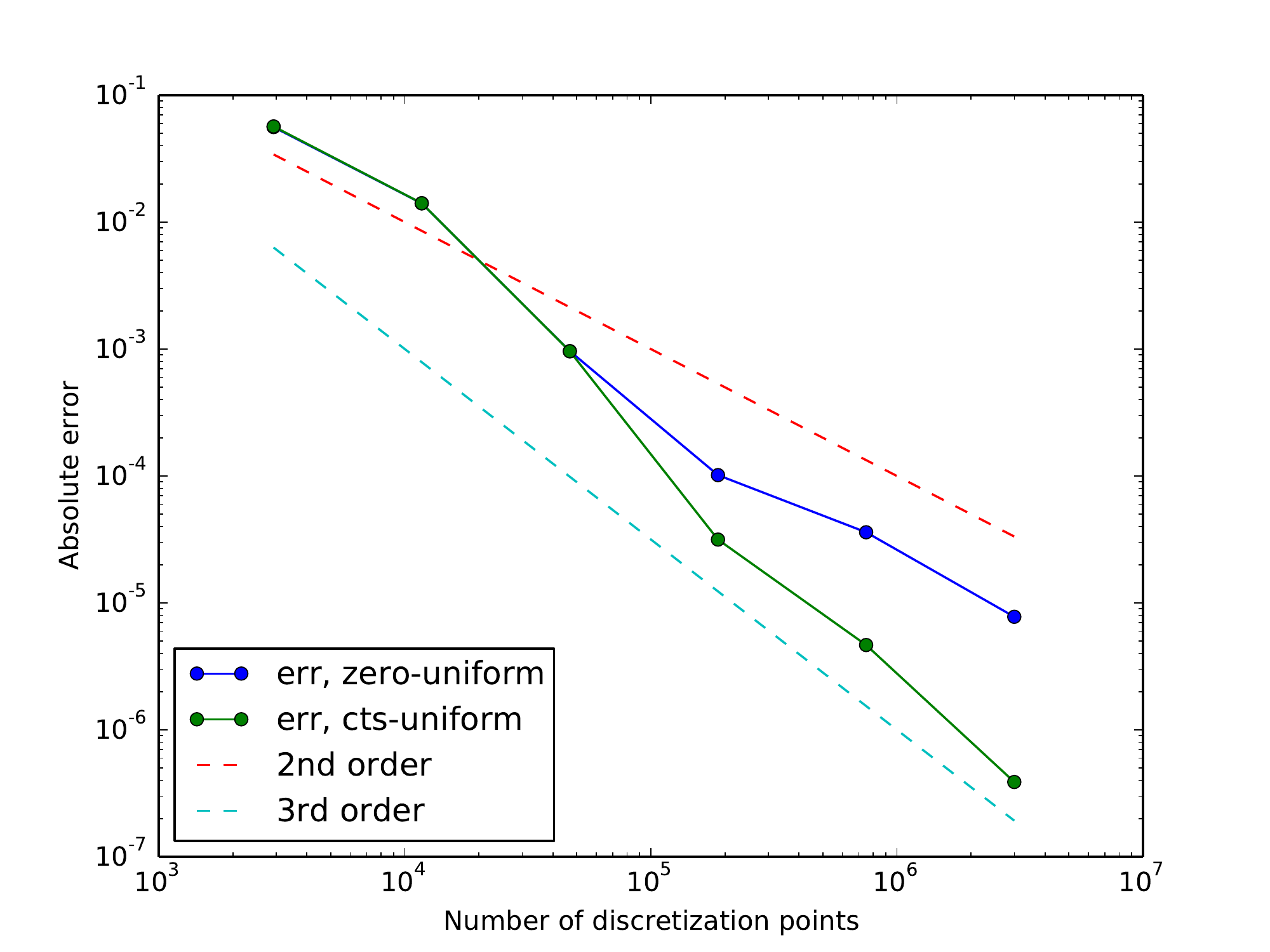}
    \includegraphics[width=.45\textwidth]{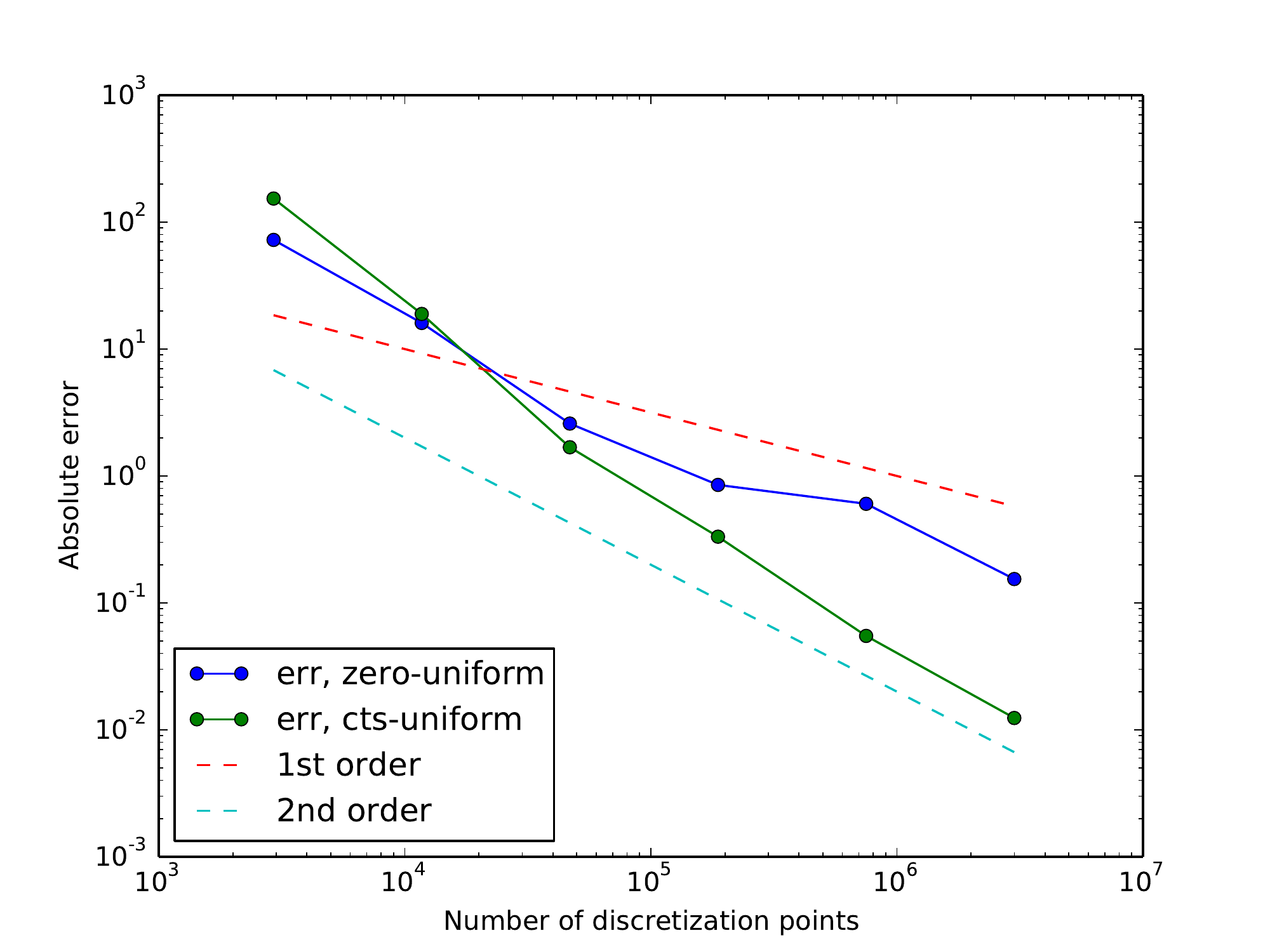}
    \caption{Example 2, convergence rates
      on a uniform tree. Accuracy of the potential (left) and gradient (right) versus the number of discretization nodes $N_\Omega$, using either extension-by-zero (blue curves) or continuous extension (green curves).}
  \label{fig:ex2uni}
  \end{centering}
\end{figure}

Next, we consider the question of the convergence order using extension-by-zero and continuous extension with a  layer potential. In Figure \ref{fig:ex2uni}, we plot the
error for the potential and gradient for increasing $N_\Omega$ on uniform trees with both extension-by-zero and continuous extension. For this example, the two methods
have similar error until $N_\Omega$ is large because the  irregularity in the solution is unresolved by the grid for small $N_\Omega$. Once $N_\Omega$ is sufficiently large, we see that the convergence rate for continuous extension is faster, though the specific rates are not as clear as they were for Example 1.

\begin{figure}[h!]
  \begin{centering}
    \includegraphics[width=.9\textwidth]{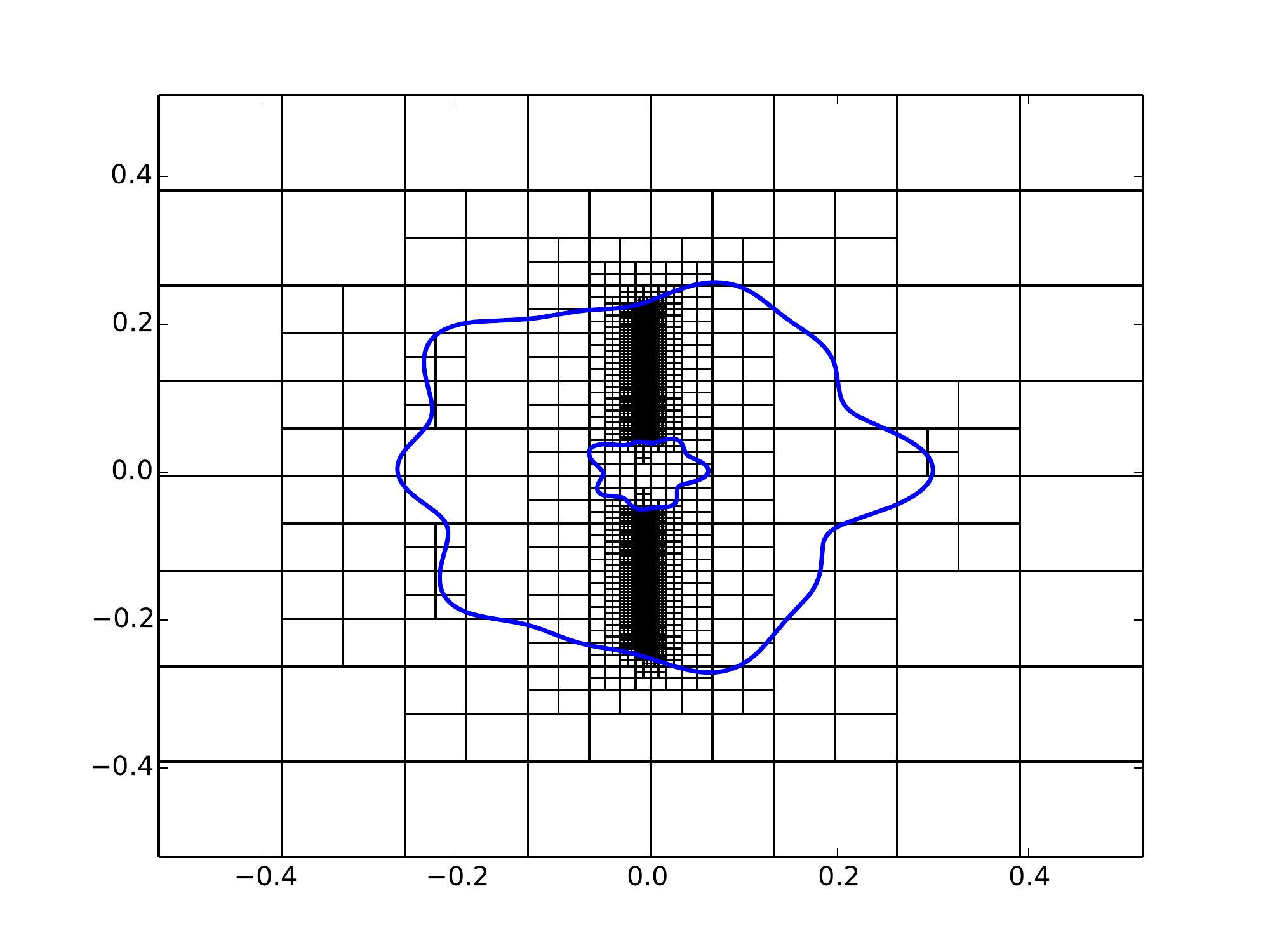}
    \caption{Example 2. An example adaptive tree for 
      the continuously extended $f_e$.}
  \label{fig:ex2adpboxes}
  \end{centering}
\end{figure}

Figure \ref{fig:ex2uni} also demonstrates that a uniform grid does a poor job of giving high accuracy for the gradient. We now test the effect of adaptive mesh refinement as in Example 1. Figure \ref{fig:ex2adpboxes} shows a representative adaptive tree for Example 2. The \textit{a priori} refinement strategy places many boxes near the irregularity in $f_e$. 
Because the continuous extension is smooth outside of $\Omega$, the effect of the ``ridge'' on the $x_2$ axis does not  extend far outside of the domain. 

\begin{figure}[h!]
  \begin{centering}
    \includegraphics[width=.45\textwidth]{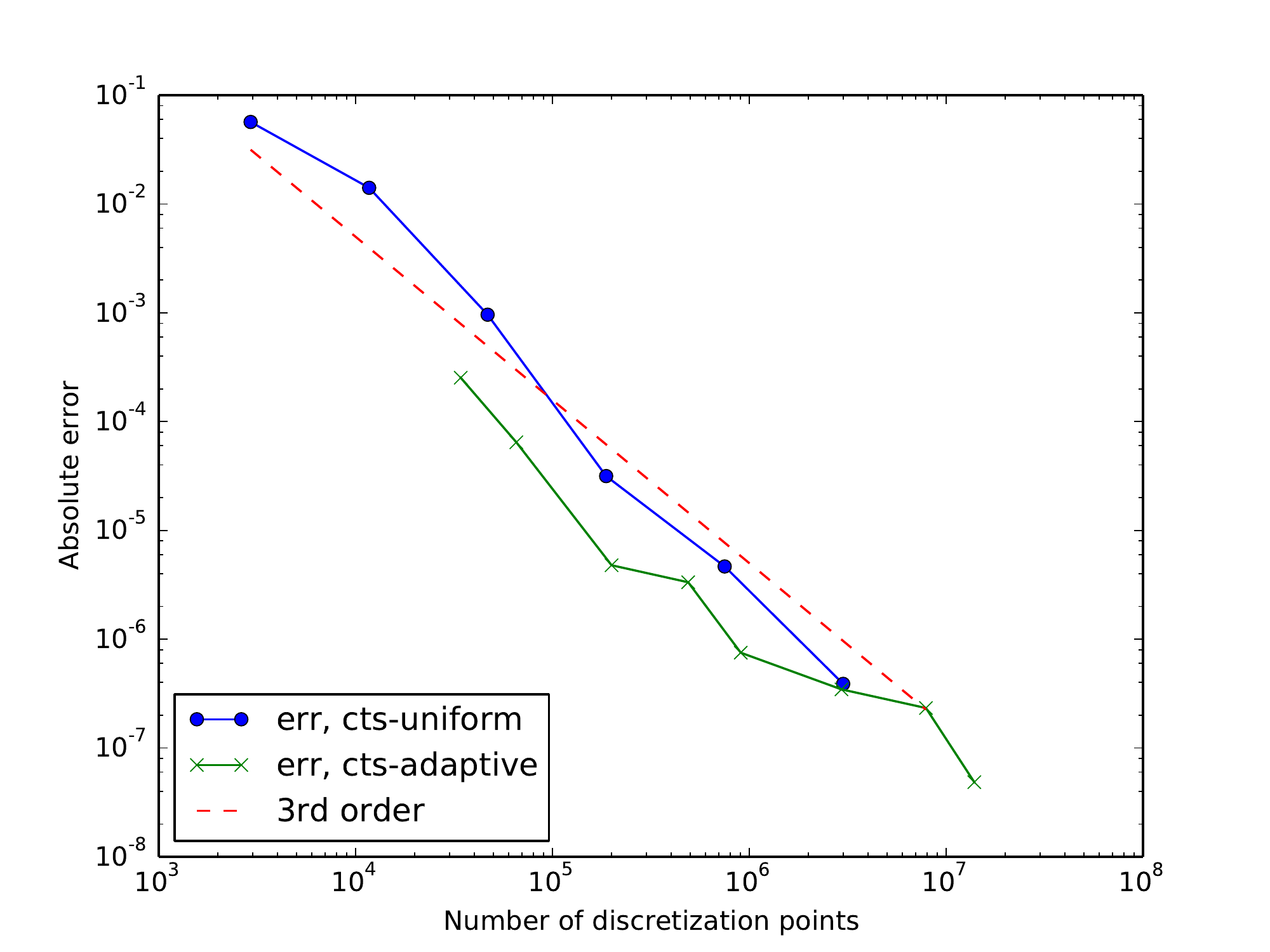}
    \includegraphics[width=.45\textwidth]{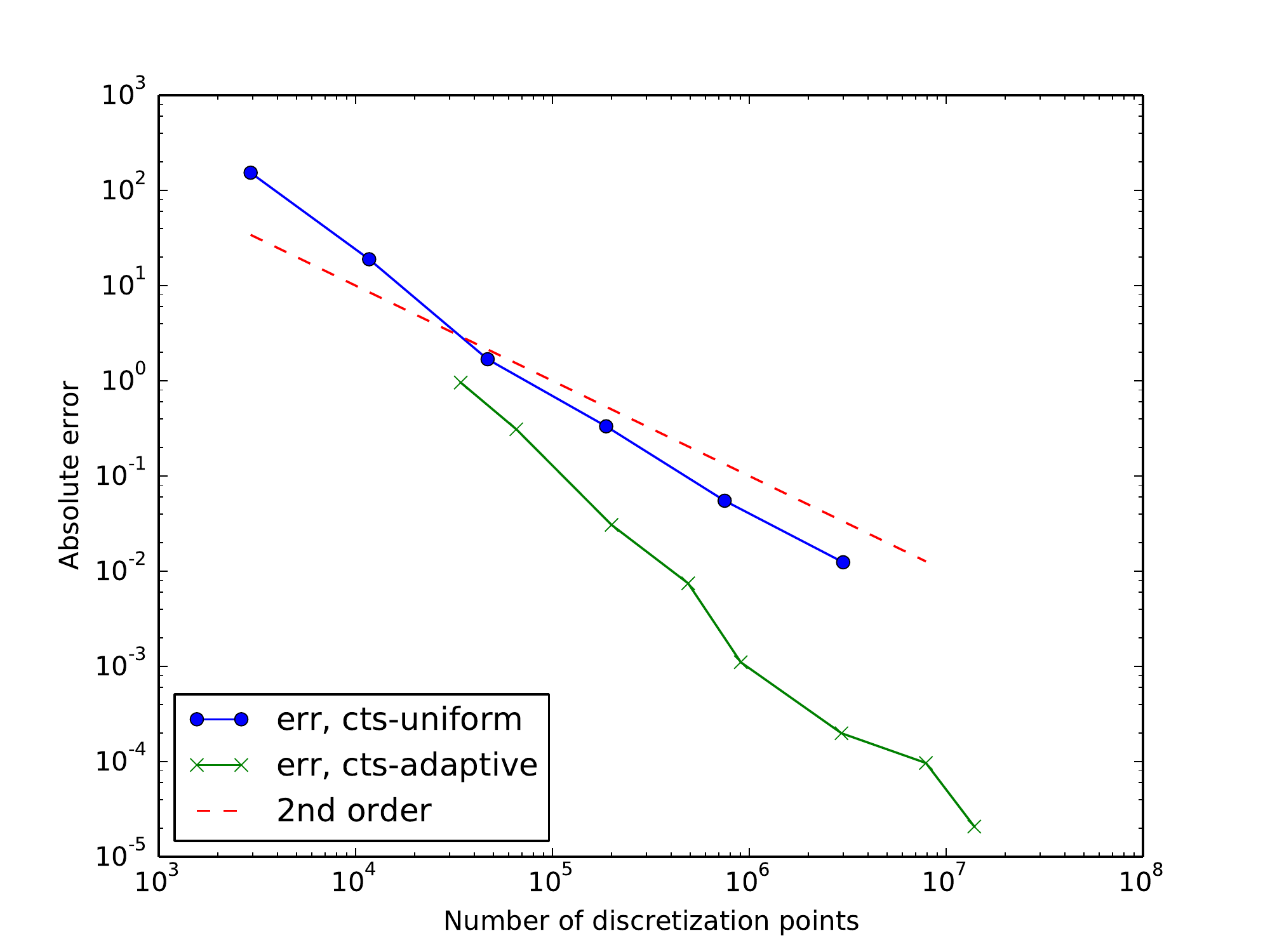}
    \caption{Example 2. Using continuous extension,  
      a plot of the error in the potential (left) and error 
      in the gradient (right) versus the number of discretization 
      nodes $N_\Omega$,
      for both uniform (blue dots) and adaptive trees (green crosses).}
  \label{fig:ex2adp}
  \end{centering}
\end{figure}

As in Example 1, the adaptive discretization strategy provides modest improvement in the accuracy of the potential (and eventually no improvement at all). An explanation for this is that the solution $u$ is much smoother than $f$, so we greatly over-resolve $u$ when we construct the tree with the goal of resolving $f$. In other words, the \textit{a priori} refinement strategy is eventually less efficient than the uniform strategy in terms of the accuracy of the potential. In contrast, adaptive discretization provides significant gains
in the accuracy of the gradient. This is because the gradient is less smooth and more difficult to resolve than $u$ so that the additional boxes used to resolve $f$ are not
as wasteful. Note also that once again, for the tests with adaptive grids much larger values of $N_\Omega$ could be achieved, for the same reasons as in Example 1.

Finally, we present run time performance results as 
we did for Example 1. The conclusions here are the same as the ones for Example 1. Observe in particular that the performance of the box code is nearly the same here as it was for Example 1, even though the trees used in this example are highly adaptive. This is one of the major advantages of the numerical method we present in this article. 

\begin{table}
\centering
\begin{tabular}{lllll} 
\toprule 
$N_\Omega$ & $N_V$ & $t_V$ & $N_\Omega/t_V$ & $N_V/t_V$ \\ 
\midrule 
2.9290e+03 & 1.6384e+04 & 1.0521e-01 & 2.7840e+04 & 1.5573e+05 \\ 
1.1717e+04 & 6.5536e+04 & 4.5237e-02 & 2.5901e+05 & 1.4487e+06 \\ 
4.6846e+04 & 2.6214e+05 & 1.1961e-01 & 3.9166e+05 & 2.1917e+06 \\ 
1.8739e+05 & 1.0486e+06 & 3.9992e-01 & 4.6856e+05 & 2.6220e+06 \\ 
7.4955e+05 & 4.1943e+06 & 1.4935e+00 & 5.0188e+05 & 2.8084e+06 \\ 
2.9983e+06 & 1.6777e+07 & 7.0171e+00 & 4.2728e+05 & 2.3909e+06 \\ 
\bottomrule 
\end{tabular}
\caption{Box code timing information for Example 2 with continuous
function extension and a uniform tree.\label{tab:ex2unitimes}}
\end{table}

\begin{table}
\centering
\begin{tabular}{lllll} 
\toprule 
$N_\Omega$ & $N_V$ & $t_V$ & $N_\Omega/t_V$ & $N_V/t_V$ \\ 
\midrule 
3.4204e+04 & 3.8032e+04 & 4.5255e-02 & 7.5581e+05 & 8.4039e+05 \\ 
6.5547e+04 & 7.0480e+04 & 5.4312e-02 & 1.2069e+06 & 1.2977e+06 \\ 
1.9972e+05 & 2.0987e+05 & 1.1204e-01 & 1.7826e+06 & 1.8732e+06 \\ 
4.8924e+05 & 5.1256e+05 & 2.3064e-01 & 2.1212e+06 & 2.2223e+06 \\ 
9.0490e+05 & 9.6006e+05 & 4.0450e-01 & 2.2371e+06 & 2.3735e+06 \\ 
2.9398e+06 & 3.0676e+06 & 1.1827e+00 & 2.4857e+06 & 2.5937e+06 \\ 
\bottomrule 
\end{tabular}
\caption{Box code timing information for Example 2 with continuous
function extension and an adaptive tree.
\label{tab:ex2adptimes}}
\end{table}

\begin{table}
\centering
\begin{tabular}{lllll} 
\toprule 
$N_\Omega$ & $t_{QP}$ & $t_{QE}$ & $N_\Omega/(t_{QP}+t_{QE})$ & $N_\Omega/t_{QE}$ \\ 
\midrule 
2.9290e+03 & 2.5967e+00 & 1.2099e-03 & 1.1274e+03 & 2.4209e+06 \\ 
1.1717e+04 & 1.8191e+00 & 3.3716e-03 & 6.4292e+03 & 3.4752e+06 \\ 
4.6846e+04 & 1.7931e+00 & 1.2859e-02 & 2.5940e+04 & 3.6431e+06 \\ 
1.8739e+05 & 1.7918e+00 & 4.8408e-02 & 1.0183e+05 & 3.8710e+06 \\ 
7.4955e+05 & 1.7931e+00 & 2.0054e-01 & 3.7597e+05 & 3.7377e+06 \\ 
2.9983e+06 & 1.8273e+00 & 7.3797e-01 & 1.1688e+06 & 4.0629e+06 \\ 
\bottomrule 
\end{tabular}
\caption{QBX timing information for Example 2 with 
continuous function extension and a uniform tree.
\label{tab:ex2uniqbxtimes}}
\end{table}

\begin{table}
\centering
\begin{tabular}{lllll} 
\toprule 
$N_\Omega$ & $t_{QP}$ & $t_{QE}$ & $N_\Omega/(t_{QP}+t_{QE})$ & $N_\Omega/t_{QE}$ \\ 
\midrule 
3.4204e+04 & 2.1972e+00 & 1.1984e-02 & 1.5483e+04 & 2.8541e+06 \\ 
6.5547e+04 & 2.1816e+00 & 2.5388e-02 & 2.9700e+04 & 2.5818e+06 \\ 
1.9972e+05 & 2.1215e+00 & 5.0410e-02 & 9.1958e+04 & 3.9620e+06 \\ 
4.8924e+05 & 2.0789e+00 & 1.2520e-01 & 2.2197e+05 & 3.9077e+06 \\ 
9.0490e+05 & 1.8968e+00 & 2.2485e-01 & 4.2651e+05 & 4.0245e+06 \\ 
2.9398e+06 & 2.1047e+00 & 7.7307e-01 & 1.0216e+06 & 3.8028e+06 \\ 
\bottomrule 
\end{tabular}
\caption{QBX timing information for Example 2 with 
continuous function extension and an adaptive tree.
\label{tab:ex2adpqbxtimes}}
\end{table}

\section{Conclusion}
\label{sec:conclusion}
We have demonstrated that continuous global function extension constructed as the solution of an exterior Laplace problem provided an effective framework to apply adaptive FMM based Poisson solvers to problems with complex geometries. We found that the desirable properties of the FMM are kept intact with such a method: the amount of work still scales linearly with the number of degrees of freedom in the computational domain and is competitive with classical FFT-based solvers in terms of work per grid point, despite the flexibility of adaptive mesh refinement. This holds even for multiply connected domains with irregular boundaries. The adaptive refinement capability of our new solver plays a crucial role in guaranteeing an efficient use of the degrees of freedom in the system, and in obtaining high accuracy for the gradient of the potential. Finally, for the particular situations in which a smooth global extension is readily available without resorting to numerical computation, as is for example the case of an extension by zero in plasma physics applications \cite{Lee}, we have presented a numerical method which leads to the same order of convergence for the gradient of the potential as the potential itself. In our implementation of the FMM, this translates to 4th order convergence for both the potential and the gradient, and the order of convergence can be increased by choosing higher order basis functions \cite{Ethridge}. 

Of course, when continuous extension is employed, the 
{\em convergence order} of the method is not 
particularly high. We demonstrated above that adaptive refinement 
can help improve the accuracy per degree of freedom in this case,
particularly for the gradient, but the low order of accuracy is 
really a result of compromise. The method of this paper emphasizes
ease of use, domain flexibility, speed, and compatibility with
adaptive refinement strategies. To achieve these goals we have chosen
an embedded boundary method (for ease of use and domain flexibility) built 
on a box code (for speed and handling highly adaptive grids).
Because it is an embedded boundary method, high order accuracy
is more difficult to achieve.
However, the method asks for very little from the user. Only a
parametric description of the boundary and a method for 
evaluating $f$ accurately in the domain must be provided. In particular,
no special quadrature rules are required, as is the case for a 
boundary fitted mesh, and there are no requirements on the accuracy
of derivatives of the user-provided $f$. As noted in Section
\ref{subsec:funextend}, when accurate derivatives
of $f$ are available, an extension computed as the solution of a 
polyharmonic equation would result in a higher order method.

The capabilities of our solver can be extended in a number of ways. First, $C^{1}$ function extension provided by the solution of an exterior biharmonic problem would lead to faster convergence for the solution and gradient than we have obtained with $C^{0}$ extension, provided that accurate values for the gradient of $f$ are available on the boundary. Second, one could allow for boundaries with corners and which nearly self intersect. Numerical tools addressing these two challenges have recently been developed, but have not yet been implemented in the Poisson context. Fortunately, the overall method is largely agnostic as to how the function extension and harmonic correction are computed, so that new methods may be swapped in when they become available. Finally, much of the technology and analysis required for this work extends to three dimensions in a straightforward manner. This is the subject of ongoing work, with progress to be reported at a later date.

\section{Acknowledgments}
The authors would like to thank Prof. Leslie Greengard (NYU) and Dr. Manas Rachh (Yale) for many insightful conversations, and Dr. Zydrunas Gimbutas (NIST) for helping with the generation of tables. T.A. was partially supported by the U.S. Department of Energy under contract DEFG0288ER25053, by the Air Force Office of Scientific Research under NSSEFF Program Award FA9550-10-1-0180 and FA9550-15-1-0385, and by a GSAS Dissertation Fellowship from NYU. A.J.C. was supported by the U.S. Department of Energy, Office of Science, Fusion Energy Sciences under Award Nos. DE-FG02-86ER53223 and DE-SC0012398.

%%%%%%%%%%%%%%%%%%%%%%%%%%%%%%%%%%%%%%%%%%%%%%%%%%%%%%%%%%%%%%%

%%%%%%%%%%%%%%%%%%%%%%%%%%%%%%%%%%%%%%%%%%%%%%%%%%%%%%%%%%%%%%%
%%%%%%%%%%%%%%%%%%%%%%%%%%%%%%%%%%%%%%%%%%%%%%%%%%%%%%%%%%%%%%%
%%%%%%%%%%%%%%%%%%%%%%%%%%%%%%%%%%%%%%%%%%%%%%%%%%%%%%%%%%%%%%%
%%%%%%%%%%%%%%%%%%%%%%%%%%%%%%%%%%%%%%%%%%%%%%%%%%%%%%%%%%%%%%%
%%%%%%%%%%%%%%%%%%%%%%%%%%%%%%%%%%%%%%%%%%%%%%%%%%%%%%%%%%%%%%%
\begin{appendix} 

\end{appendix}

\bibliographystyle{elsarticle-num}
% \bibliography{refs}

\end{document}